\documentclass[preprint]{elsarticle}

\usepackage{amssymb}
\usepackage{graphicx,amsmath,bm,color,geometry,subfig}
\usepackage[pagewise]{lineno}
\usepackage{algorithmic,algorithm,diagbox} 
\usepackage{tabularx}
\usepackage{booktabs}
\usepackage{multirow}

\newcommand{\norv}{m}               
\newcommand{\pvc}{x}
\newcommand{\svc}{\xi}
\newcommand{\pv}{\bm{\pvc}}
\newcommand{\sv}{\bm{\svc}}

\newcommand{\rd}{\mathrm{d}}

\newcommand{\trans}{{\scriptsize\mbox{T}}}
\newcommand{\tol}{\mathtt{TOL}}

\newcommand{\ffemsoln}{u_h}                                   	
\newcommand{\vfemsoln}{\bm{u}_{\sv}}                          	
\newcommand{\frbfemsoln}{u^{(r)}_h}                           	
\newcommand{\vrbfemsoln}{\bm{u}^{(r)}_{\sv}}                	
\newcommand{\fsfemsoln}{u_{hp}}                               	

\newcommand{\frbsfemsoln}{u^{(r)}_{hp}}                     	
\newcommand{\mrbsfemsoln}{\bm{U}^{(r)}}                        	
\newcommand{\vrbsfemsoln}{\bm{u}^{(r)}}                        	
\newcommand{\frefsoln}{u_h^{\scriptsize\mbox{ref}}}           	

\newcommand{\ofrefsoln}{u_{\scriptsize\mbox{ref}}}           	

\newcommand{\pE}[1]{\mathbb{E}[#1]}
\newcommand{\pV}[1]{\mathbb{V}[#1]}
\newcommand{\meanerr}{\mathtt{err}_{\mbox{m}}}
\newcommand{\varerr}{\mathtt{err}_{\mbox{v}}}

\newcommand{\Nx}{N_{h}}

\newcommand{\Nxi}{N_{p}}

\newcommand{\mypdf}{\rho}
\newcommand{\smp}{\bm{\varXi}}

\makeatletter
\def\ps@pprintTitle{%
	\let\@oddhead\@empty
	\let\@evenhead\@empty
	\let\@oddfoot\@empty
	\let\@evenfoot\@oddfoot
}
\makeatother

\journal{Applied Mathematics and Computation}

\begin{document}

\begin{frontmatter}

\title{Reduced basis stochastic Galerkin methods for partial differential equations with random inputs}
\author[label1]{Guanjie Wang}
\ead{guanjie@lixin.edu.cn}
\affiliation[label1]{organization={School of Statistics and Mathematics, Shanghai Lixin University of Accounting and Finance},
	city={Shanghai},
	postcode={201209},
	country={P. R. China}}

\author[label2]{Qifeng Liao\corref{qfcor}}
\ead{liaoqf@shanghaitech.edu.cn}
\cortext[qfcor]{Corresponding author.}
\affiliation[label2]{organization={School of Information Science and Technology, ShanghaiTech University},
	city={Shanghai},
	postcode={201210},
	country={P. R. China}}

\begin{abstract}
We present a reduced basis stochastic Galerkin method for partial differential equations with random inputs. In this method, the reduced basis methodology is integrated into the stochastic Galerkin method, resulting in a significant reduction in the cost of solving the Galerkin system. To reduce the main cost of matrix-vector manipulation involved in our reduced basis stochastic Galerkin approach, the secant method is applied to identify the number of reduced basis functions. We present a general mathematical framework of the methodology, validate its accuracy and demonstrate its efficiency with numerical experiments.
\end{abstract}




\end{frontmatter}


\section{Introduction}
In the past decades, there has been rapid development in efficient numerical methods for solving partial differential equations (PDEs) with random inputs. This explosion of interest has been driven by the need to conduct uncertainty quantification for modeling realistic problems, such as diffusion problems and acoustic scattering problems. The sources of uncertainty for these problems typically arise from a lack of knowledge or measurements of realistic model parameters, such as permeability coefficients and refraction coefficients.

It is of great interest to compute the mean and variance of solutions for PDEs with random inputs. To this end, a lot of efforts are made. One of the simplest ways to deal with this problem is the Monte Carlo method (MCM) and its variants~\cite{Caflisch1998,fishman2013}. In MCM, large number of realizations of the random inputs are generated  based on their probability density functions. For each realization, the associated deterministic problem can be solved by available numerical schemes. The mean and variance of the stochastic solution can then be approached  by the statistical information obtained from these solutions to the deterministic problems. While the MCM is straightforward to implement, the convergence of solution statistics is relatively slow, typically requiring numerous realizations.

To improve the efficiency of sampling based methods, the stochastic collocation method (SCM) based on sparse grids is developed in \cite{Xiu05} (and a comprehensive review can be found in \cite{Xiu2010}). After that, the reduced basis collocation method (RBCM) is proposed in \cite{Elmanliao}, which significantly reduces the computational costs associated with collocation methods without loss of accuracy. Meanwhile, reduced basis methods are actively developed for solving PDEs with random inputs \cite{Boyaval2010Reduced,Quarteroni2016Reduced,Chen2019resfree,Chen2014Comparison}. Notably, efficient preconditioning techniques for reduced basis methods are developed in \cite{Elman2015pcond}.

An alternative efficient approach for handling PDEs with random inputs is the stochastic Galerkin method (SGM)~\cite{Ghanem2003, Xiu2002wiener,Xiu2002modeling}. In stochastic Galerkin methods, the stochastic solution is represented by a finite expansion of trial basis functions, and the expansion coefficients can be obtained by solving a linear system resulting from the Galerkin projection. The trial functions primarily include polynomial chaos (PC)~\cite{Ghanem2003}, generalized polynomial chaos (gPC)~\cite{Xiu2002modeling,Xiu2002wiener}, piecewise polynomial bases~\cite{Babuska2004}, multi-element generalized polynomial chaos (ME-gPC)~\cite{Wan2005adaptive,Wan2006multi}, and dynamically bi-orthogonal polynomials~\cite{Cheng2013a,Cheng2013b,Musharbash,Zhoutaoa}. In this work, we mainly focus on the gPC approach for the stochastic solution. 

Designing efficient solvers is a crucial and challenging problem when dealing with the stochastic Galerkin method, as it often leads to a large coupled linear system. Various iterative solvers such as the mean based preconditioning iterative method \cite{Powell2009}, the reduced basis solver based on low-rank approximation~\cite{powell2015}, and the low-rank iterative methods~\cite{LeeElman16}, have been extensively studied. In this work, we develop a reduced basis stochastic Galerkin method (RBSGM) that integrates the reduced basis methodology into the stochastic Galerkin method, resulting in a significant reduction in the cost of solving the Galerkin system. To further minimize the primary computational cost in our RBSGM, which encompasses matrix-vector manipulations throughout the iterative procedure, we employ the secant method to ascertain the optimal number of reduced basis functions. While the reduced basis method for physical approximation in stochastic Galerkin methods has been previously explored in \cite{fluids6080263}, our contribution lies in presenting a systematic procedure that integrates the reduced basis methodology into the stochastic Galerkin method. This approach leads to a substantial reduction in the computational cost associated with matrix-vector multiplications.

The outline of this work is as follows: In Section 2, we describe the PDE formulation considered in this study, introduce the stochastic Galerkin method, and discuss the iterative method used for solving the stochastic Galerkin system. In Section 3, we present our reduced basis stochastic Galerkin method. Section 4 provides the numerical results, and Section 5 concludes the paper.

\section{Problem setting and stochastic Galerkin method}\label{sec:sg}
This section describes the PDEs with random inputs considered in this study, introduces the stochastic Galerkin method, and discusses the iterative method used for solving the stochastic Galerkin system.

\subsection{Problem setting}
Let $\sv = \left[\svc_{1},\ldots,\svc_{\norv}\right]^{\trans}$ be an $\norv$-variable random vector. The image of $\svc_i$ is denoted by $\Gamma_i$, and the probability density function of $\svc_i$ is given by $\rho_i(\svc_i)$. If we further assume that the components of $\sv$, i.e., $\svc_{1},\ldots,\svc_{\norv}$, are mutually independent, then the image of $\sv$ is given by $\Gamma = \Gamma_1\times\cdots\times\Gamma_{\norv}$, and the probability density function of $\sv$ is given by $\mypdf(\sv) = \prod_{i=1}^{\norv}\mypdf_i(\svc_i)$. %

In this work, we consider the following partial differential equations (PDEs) with random inputs, which are widely used in modeling steady state diffusion problems~\cite{Xiu2002modeling} and  acoustic scattering problems~\cite{Elman2005,xiushen07}
\begin{equation}\label{eq:spde}
\begin{cases}
-\nabla\cdot [a(\pv,\sv)\nabla u(\pv,\sv)]-\kappa^2(\pv,\sv)u(\pv,\sv)=f(\pv)\quad & \forall (\pv,\sv)\in D\times\Gamma,\\
\mathfrak{b}(\pv,\sv,u(\pv,\sv))=g(\pv)\quad &\forall (\pv,\sv)\in\partial D\times \Gamma.
\end{cases}
\end{equation}
Here $\mathfrak{b}$ is a boundary operator, $f$ is the source function and $g$ specifies the boundary condition. Both  $a(\pv,\sv)$ and $\kappa(\pv,\sv)$ are assumed to have the following form:
\begin{equation}\notag
a(\pv,\sv) =  \sum_{i=0}^{m}a_0(\pv)\svc_{i},  \ 
\kappa(\pv,\sv) =  \sum_{i=0}^{m}\kappa_i(\pv)\svc_{i},\  \xi_0 = 1.
\end{equation}
Note that we denote $1$ as $\xi_0$ for simplicity. The  problem \eqref{eq:spde} is a diffusion equation when $\kappa(\pv,\sv)\equiv 0$ and  a Helmholtz equation when $a(\pv,\sv)\equiv 1$.

To simplify the presentation, we assume that the problem \eqref{eq:spde} satisfies the homogeneous Dirichlet boundary condition. However, it should be noted that the approach can be straightforwardly extended to handle non-homogeneous boundary conditions.

\subsection{Variational formulation}
To introduce the  variational form of \eqref{eq:spde}, some notations are required.
We first consider the Hilbert space 
\begin{equation}\notag\label{eq:l2space}
L^2(D): = \left\{v: D \to \mathbb{C}\ \bigg|\ \Vert v\Vert_{L^2} <\infty \right\}, 
\end{equation}
of square integrable functions equipped with the inner product 
\begin{equation}\notag
\langle u, v\rangle_{L^2} = \int_{D}v^*u\,\rd \pv,
\end{equation}
and the norm
\begin{equation}\notag\label{eq:l2norm}
\Vert v \Vert_{L^2}:=\sqrt{\langle v,v\rangle_{L^2}}.
\end{equation}
Moreover, let 
\begin{equation}\notag
H_0^1(D):=\left\{v\in H^1(D)\,| \,  v=0 \textrm{ on } \partial D_D\right\},
\end{equation}
where $H^1(D)$ is the Sobolev  space
\begin{equation}\notag
H^1(D):=\left\{v\in L^2(D)\,, \,   \partial v/ \partial x_i\in L^2(D),  i=1,\ldots,d\right\}.
\end{equation}
Next, we define the Hilbert space
\begin{equation}\notag
L_{\mypdf}^2(\Gamma):=\left\{v(\sv): \Gamma \to \mathbb{R}\ \bigg|\ \Vert v\Vert^2_{L^2_{\mypdf}} <\infty \right\},
\end{equation}
which is equipped with the inner product
\begin{equation}\notag\label{eq:inner_stoch}
\langle u,v\rangle_{L^2_\mypdf}=\int_{\Gamma}\mypdf(\sv)v^*(\sv)u(\sv)\,\rd\sv,
\end{equation}
and the norm 
\begin{equation}\notag
\Vert v \Vert_{L^2_{\mypdf}}:=\sqrt{\langle v,v\rangle_{L^2_{\mypdf}}}.
\end{equation}
Following \cite{Babuska2004}, the tensor product space of $L^2(D)$ and $L^2_{\mypdf}(\Gamma)$ is defined as: 
\begin{equation}\notag
	L^2(D)\otimes L^2_{\mypdf}(\Gamma) :=\left\{ w(\pv,\sv) \bigg| w(\pv,\sv)=\sum_{i=1}^{n}u_i(\pv)v_i(\sv), u_i(\pv)\in L^2(D),\, v_i(\sv)\in L^2_{\mypdf}(\Gamma), n\in \mathbb{N}^{+} \right\},
\end{equation}
which is equipped with the inner product
\begin{equation}\notag
\left\langle w_1,w_2\right\rangle_{L^2\otimes L^2_{\mypdf}} = \int_{\Gamma}\int_{D}\mypdf(\sv) w_2^* w_1\rd \pv\,\rd \sv.
\end{equation}
Then, the variational form of \eqref{eq:spde} with respect to the inner product $\langle\cdot,\cdot\rangle_{L^2\otimes L^2_{\mypdf}}$ can be written as:
\begin{equation}\label{eq:weak}
\int_{\Gamma}\int_{D}\mypdf(\sv)\left[  a \nabla u\cdot \nabla w^*  - \kappa^2 u w^*\right]\rd \pv\rd \sv = \int_{\Gamma}\int_D \mypdf(\sv)fw^*\rd \pv\rd\sv,\ \forall w\in L^2(D)\otimes L^2_{\mypdf}(\Gamma).
\end{equation}

\subsection{Discretization}
To obtain the  discrete version of (\ref{eq:weak}), we need to introduce a finite dimensional subspace of $L^2(D)\otimes L^2_{\mypdf}(\Gamma)$ and find an approximation lies within it. Given subspace of the physical and the stochastic spaces 
\begin{equation}\notag \label{eq:sigle-bases}
V_h = {\rm span}\left\{v_s(\pv)\right\}_{s=1}^{\Nx}\subset H_0^1(D),\ 
S_p= {\rm span}\left\{\Phi_j(\sv)\right\}_{j=1}^{\Nxi}\subset L^2_{\mypdf}(\Gamma),\
\end{equation}
where $h$ is the mesh size of  the physical space and $p$ is the order of generalized polynomial chaos (gPC) of the stochastic space, the finite dimensional subspace of $L^2(D)\otimes L^2_{\mypdf}(\Gamma)$ can be defined as:
\begin{equation}\notag
W_{hp}:= V_h\otimes S_p = \mbox{span}\left\{ v_s(\pv)\Phi_j(\sv)\ \big|\ s = 1,\ldots,\Nx \mbox{ and } j = 1,\ldots,\Nxi\right\}.
\end{equation} 
In this work, the basis functions $\{v_s(\pv)\}_{s=1}^{\Nx}$ are taken to
be the standard trial (test) functions of finite element approximations. Meanwhile, $\left\{\Phi_j(\sv)\right\}_{j=1}^{\Nxi}$ are taken to be the generalized polynomial chaos (gPC) basis, which consists of the orthonormal polynomials with respect to the inner product~$\langle \cdot,\cdot\rangle_{L^2_\mypdf}$, i.e.,
\begin{equation}\notag
\left\langle\Phi_j(\sv),\Phi_k(\sv)\right\rangle_{L^2_{\mypdf}}=\int_{\Gamma} \mypdf(\sv)\Phi_j(\sv)\Phi_k^*(\sv)\, \rd\sv= \delta_{jk},
\end{equation}
where $\delta_{jk}$ is the Kronecker delta function. 
As usual, $\left\{\Phi_j(\sv)\right\}_{j=1}^{\Nxi}$
consists of orthonormal polynomials with
total degrees up to $p$, and $p$ is referred to as the gPC order. Then the dimension of $S_p$ is given by $\Nxi=(m+p)!/(m!p!)$, , where $m$ is the number of random variables. For more details about gPC, the readers can consult Refs.~\cite{Xiu2010,Xiu2002wiener}.

Suppose that $\fsfemsoln(\pv,\sv)$ is an approximation of $u(\pv,\sv)$  that  belongs to $W_{hp}$, i.e., 
\begin{equation}\label{eq:uapp}
\fsfemsoln(\pv,\sv) = \sum_{j=1}^{\Nxi}\sum_{s=1}^{\Nx} u_{js}v_s(\pv)\Phi_j(\sv).
\end{equation}
By restricting the test functions in the weak form \eqref{eq:weak} to $W_{hp}$, we have 
\begin{equation}\notag
\int_{\Gamma}\int_{D}\mypdf(\sv)\left[  a \nabla \fsfemsoln\cdot \nabla w^*  - \kappa^2 \fsfemsoln w^*\right]\rd \pv\rd \sv = \int_{\Gamma}\int_D \mypdf(\sv)fw^*\rd \pv\rd\sv,\ \forall  w\in W_{hp}.
\end{equation}
This gives rise to a linear system
\begin{equation}\label{eq:linear_sg}
\bm{A}{\bm{u}} = \bm{b},
\end{equation}
where ${\bm{u}}^\trans = [\bm{u}_1^\trans,\ldots,\bm{u}_{\Nxi}^\trans]$, $\bm{u}_j\in R^{\Nx\times 1}$, $\bm{u}_j(s) = u_{js}$ and 
\begin{equation}\label{eq:linear_lhs}
\bm{A} = \sum_{i=0}^{m}\bm{G}_{i0}\otimes \bm{A}_i -\sum_{i=0}^{m}\sum_{j=0}^{m}\bm{G}_{ij} \otimes \bm{B}_{ij}, \ \bm{b}=\bm{h}\otimes\bm{f}.
\end{equation}
In \eqref{eq:linear_lhs}, $\bm{G}_{ij}\in R^{\Nxi\times \Nxi}$ and $\bm{h}\in R^{\Nxi\times 1}$ only depend on the basis functions in the stochastic space. Therefore, they are referred to as the stochastic Galerkin matrices and vector. They are given by
\begin{equation}\notag\label{eq:matrix_stoch}
\bm{G}_{ij}(l,n) = \int_{\Gamma}\mypdf(\sv)\svc_i\Phi_l(\sv)\svc_j\Phi_n^*(\sv)\rd\sv, \ \bm{h}(l) = \int_{\Gamma}\mypdf(\sv){\Phi_l(\sv)}\rd\sv,\ \xi_0=1,
\end{equation}
where $i,j\in \{0,1,\ldots,m\}$ and $l,n\in \{1,2,\ldots,\Nxi\}$. 

On the other hand, the matrices $\bm{A}_i$, $\bm{B}_{ij}$ and the vector $\bm{f}$  are given by
\begin{eqnarray}\label{eq:matrix_phy}
\bm{A}_i(s,t) = \int_{D} a_i(\pv)\nabla v_s(\pv)\cdot\nabla v_t^*(\pv)\rd \pv,\ \bm{B}_{ij}(s,t) = \int_{D}\kappa_i\kappa_jv_sv_t^*\rd \pv,\ \bm{f}(s) = \int_{D} fv_s^*\rd \pv,
\end{eqnarray}
where $i,j\in \{0,1,\ldots,m\}$ and $s,t\in \{1,2,\ldots,\Nx\}$. It is clear that $\bm{A}_i$, $\bm{B}_{ij}$ and $\bm{f}$ only depend on the basis functions in the physical space and the coefficients $a_i(\pv)$, $\kappa_i(\pv)$.

Once the unknowns are solved through \eqref{eq:linear_sg}, the approximation $\fsfemsoln(\pv,\sv)$ can be easily reconstructed by \eqref{eq:uapp}. Since $\{\Phi_j(\sv)\}_{j=1}^{\Nxi}$ is orthonormal, the mean and  the variance of $\fsfemsoln(\pv,\sv)$ are given by
\begin{equation}\notag
\pE{\fsfemsoln} = \sum_{s=1}^{\Nx}u_{1s}v_s(\pv),\ \pV{\fsfemsoln} = \sum_{j=2}^{\Nxi}\sum_{s=1}^{\Nx} (u_{js}v_s(\pv))^2.
\end{equation}

\subsection{Iterative methods}
The coefficient matrix $\bm{A}$  in  \eqref{eq:linear_sg} is usually a block-wise sparse matrix~\cite{Ernst2010,Liao2019}. Note that the size of matrix $\bm{A}$ in \eqref{eq:linear_sg} is $\Nx \Nxi$, where $\Nx$ is the number of degrees of freedom (DOF) used in the spatial discretization, and $\Nxi$ is the number of basis used in the stochastic space. To achieve a high accuracy approximation, both $\Nx$ and $\Nxi$ are often very large, resulting in a large sparse linear system that needs to be solved. In this section, we discuss  iterative methods for this kind of linear systems. In particular, we focus on Krylov subspace methods~\cite{Golub2013,Elman2014}, which are based on the projection of the linear system  (\ref{eq:linear_sg}) into a consecutively constructed Krylov subspace
\begin{equation}\notag
\mathcal{K}_n(\bm{A},\bm{r}^{(0)})={\rm span}\{\bm{r}^{(0)},\bm{A}\bm{r}^{(0)},\bm{A}^2\bm{r}^{(0)},\ldots,\bm{A}^{n-1}\bm{r}^{(0)}\},
\end{equation}
where
\begin{equation}\notag 
\label{eq:RES}
\bm{r}^{(0)}=\bm{b}-\bm{A}{\bm{u}}^{(0)}.
\end{equation}
and ${\bm{u}}^{(0)}$ is an initial approximation.

In Krylov subspace methods, we usually need to compute the matrix-vector products of the form $\bm{A}\bm{x}$ and/or $\bm{A}^\trans \bm{x}$. Since the coefficient matrix has a Kronecker products structure, the matrix-vector products can be computed efficiently without explicitly assembling the coefficient matrix. To demonstrate this, we introduce the $\mathtt{vec}$ operation~\cite{Golub2013}. If~$\bm{X}\in\mathbb{R}^{m\times n}$,  then $\mathtt{vec}(\bm{X})$ is an $mn$-by-1 vector obtained by ``stacking" the columns of $\bm{X}$, i.e.,
\begin{equation}\notag
\mathtt{vec}(\bm{X}) = \begin{bmatrix}
\bm{X}(:,1) \\
\vdots\\
\bm{X}(:,n)
\end{bmatrix},
\end{equation}
where $\bm{X}(:,j)$, for $j= 1,\ldots,n$, represents the $j$-th column of $\bm{X}$. Then the matrix-vector product in the iterative methods can be rewritten as:
\begin{equation}\label{eq:mv2mm}
\bm{A}\bm{x}=\left(\sum_{i=0}^{m}\bm{G}_{i0}\otimes \bm{A}_i -\sum_{i=0}^{m}\sum_{j=0}^{m}\bm{G}_{ij} \otimes \bm{B}_{ij}\right)\bm{x} = \mathtt{vec}\left(\sum_{i=0}^{m}\bm{A}_i\bm{X}\bm{G}_{i0}^\trans  -\sum_{i=0}^{m}\sum_{j=0}^{m}\bm{B}_{ij}\bm{X}\bm{G}_{ij}^\trans\right), 
\end{equation}
where $ \mathtt{vec}(\bm{X})=\bm{x}$ and  $\bm{X}\in \mathbb{R}^{\Nx\times\Nxi}$.

Preconditioning is another key ingredient for the success of iterative methods in solving the linear system~\eqref{eq:linear_sg}. Instead of solving the original problem~\eqref{eq:linear_sg}, we solve the 
preconditioned problem
\begin{equation}\notag
\bm{A}\bm{P}^{-1}\bm{y} = \bm{b}, \mbox{ with } \bm{u} = \bm{P}^{-1}\bm{y},
\end{equation}
where the nonsingular matrix $\bm{P}$ is called a preconditioner. A lot of preconditioners that exploit the structure of the linear system are developed (see, for example \cite{Powell2009,Pellissetti2000,Sousedik2014hierarchical}). In particular, the mean based preconditioner is popular due to  its non-intrusive nature. Denote the mean value of~$\sv$ by $\sv^{(0)}$, and the identity matrix by $\bm{I}$, the mean based preconditioner is constructed based on the deterministic version of problem \eqref{eq:spde} associated with the realization $\sv^{(0)}$. Let $\bm{A}_{\sv^{(0)}}$ denote the discrete version of the differential operator associated with $\sv^{(0)}$ in \eqref{eq:spde}, then the mean based preconditioner is given by $\bm{P}=\bm{I}\otimes \bm{A}_{\sv^{(0)}}$. In this work, the mean based preconditioner is used in the iterative solvers to improve the efficiency.

The main criterion for selecting an iterative method is the knowledge of the coefficient matrix properties. For example, if $\bm{A}$ is symmetric positive definite, the conjugate gradient (CG) method is the best choice; otherwise, we should use generalized minimal residual (GMRES) method or other iterative methods. For more about this topic, the readers can consult Refs.~\cite{Barrett1994,Saad2003}.

\section{The reduced basis stochastic Galerkin method}\label{sec:RBM}
In Section \ref{sec:sg}, we introduce the framework of the stochastic Galerkin method. As we have mentioned, the stochastic Galerkin method will lead to a linear system with coefficient matrix represented as
Kronecker products. The number of unknowns is $\Nx\times\Nxi$, where $\Nx$ is the DOF in the physical space, and $\Nxi$ is the DOF in the stochastic space. 

In practical implementation, $\Nx$ is large if the spatial mesh is sufficiently fine so that the error in the physical space is acceptable. On the other hand, since $\Nxi=(m+p)!/(m!p!)$, the DOF of the stochastic space is also large if the dimension of random variables is high or high accuracy  is needed in the stochastic space. This will result in an unacceptably high cost for solving the $\Nx^2\times\Nxi^2$ linear system~\eqref{eq:linear_sg}. In this section, we show that the costs can be significantly reduced by integrating the reduced basis methodology into the iterative method.

\subsection{Reduced basis method}
The variational formulation of the  deterministic version of \eqref{eq:spde}, associated with a given $\sv$, is defined as:
\begin{equation}\notag 
\int_D\left[a\nabla u \cdot \nabla v^*-k^2uv^* \right]\rd \pv = \int_{D} f v^*\rd \pv, \ \forall v \in H_0^1(D).
\end{equation}
Let $V_h=\mathrm{span}\{v_i(\pv)\}_{i=1}^{\Nx}$ be a spatial finite element approximation space (e.g., piecewise linear or quadratic polynomial space) of dimension $\Nx$. The finite element method seeks an approximation $\ffemsoln(\pv,\sv)\in V_h$ such that  
\begin{equation}\label{eq:phy_weak3}
\int_D\left[a\nabla \ffemsoln \cdot \nabla v^*-k^2\ffemsoln v^* \right]\rd \pv = \int_{D} f v^*\rd \pv, \ \forall v \in V_h.
\end{equation}
Suppose that  
\begin{equation}\notag
\ffemsoln(\pv,\sv) = \sum_{s=1}^{\Nx}u_s(\sv)v_s(\pv),
\end{equation}
then we have the following linear system 
\begin{equation}\label{eq:linear_full}
\bm{A}_{\sv} \vfemsoln = \bm{f},
\end{equation}
where 
\begin{equation}\label{eq:femlrhs}
\bm{A}_{\sv} = \sum_{i=0}^{m}\bm{A}_i\svc_i -\sum_{i=0}^{m}\sum_{j=0}^{m}\bm{B}_{ij}\svc_i\svc_j.
\end{equation}
Note that we define $\svc_0=1$ for convenience, and $\bm{A}_i,$ $\bm{B}_{ij}$, $i,j\in \{0,\ldots,m\}$ are given by \eqref{eq:matrix_phy}.

For the reduced basis method (RBM), we first find an $r$-dimensional reduced space $Q_r \subset V_h$, where $r\ll\Nx$. We then seek a reduced solution $\frbfemsoln(\pv,\sv)\in Q_r$ that satisfies the following equation: 
\begin{equation}\label{eq:phy_rbm}
\int_D\left[a\nabla \frbfemsoln \cdot \nabla v^*-k^2\frbfemsoln v^* \right]\rd \pv = \int_{D} f v^*\rd \pv, \  \forall v\in Q_r.
\end{equation} 
Typically, the reduced problem \eqref{eq:phy_rbm} is  much smaller than the full problem \eqref{eq:phy_weak3}, allowing more efficient computation for the linear system.

Since the stochastic collocation method only requires solving the deterministic problem associated with each realization, employing RBM can significantly enhance efficiency. Furthermore, it is important to note that in the stochastic Galerkin method, the discretization of the physical space is independent of the discretization of the stochastic space. Therefore, RBM can also be applied to the stochastic Galerkin method to improve efficiency. In the following, we provide detailed information on the construction of the reduced space $Q_r$.

Assuming that a training set $\smp\subset\Gamma$ is given, such that $\{u(\pv,\sv):\sv \in \smp\}$ is accurately approximated by $\{\ffemsoln(\pv,\sv):\sv \in \smp\}$, where $\ffemsoln(\pv,\sv)$ is referred to as a snapshot. The reduced space $Q_r$ is then constructed as the span of these snapshots, given by
\begin{equation}\notag
Q_r = \mathrm{span}\{\ffemsoln(\pv,\sv^{(1)}),\ldots,\ffemsoln(\pv,\sv^{(r)})\},
\end{equation}
where $\{\sv^{(1)},\ldots,\sv^{(r)}\}$ are chosen from the training set $\smp$ using a greedy algorithm.

Let $\{q_1(\pv),\ldots,q_{r}(\pv)\}$ be the  basis of  $Q_r$, and $\bm{q}_i$ be the vector of coefficient values associated with~$q_i(\pv)$, i.e.,
\begin{equation}\label{eq:reducedbasis}
q_i(\pv) = \bm{q}_i(1)v_1(\pv) +\cdots + \bm{q}_i(\Nx)v_{\Nx}(\pv), \ i = 1,2,\ldots,r.
\end{equation}
Then the linear system for the reduced problem associated with $\sv$ is given by  
\begin{equation}\label{eq:linear_rb}
\bm{Q}_r^\trans \bm{A}_{\sv} \bm{Q}_r\vrbfemsoln = \bm{Q}_r^\trans\bm{f},
\end{equation} 
where $\bm{A}_{\sv}$ and $\bm{f}$ are given by \eqref{eq:femlrhs}, $\vrbfemsoln$ is the coefficient vector of the basis $\{q_1(\pv),\ldots,q_{r}(\pv)\}$, and
\begin{equation}\notag
\bm{Q}_r = [\bm{q}_1,\ldots,\bm{q}_{r}].
\end{equation}
By \eqref{eq:femlrhs}, \eqref{eq:linear_rb} can be rewritten as:
\begin{equation}\notag
\left(\sum_{i=0}^{m}\bm{A}_i^{(r)}\svc_i -\sum_{i=0}^{m}\sum_{j=0}^{m}\bm{B}_{ij}^{(r)}\svc_i\svc_j\right)\vrbfemsoln =\bm{f}^{(r)}.
\end{equation}
where
\begin{equation}\notag
\bm{A}_{i}^{(r)} = \bm{Q}_r^\trans\bm{A}_i\bm{Q}_r, \ \bm{B}_{ij}^{(r)} = \bm{Q}_r^\trans\bm{B}_{ij}\bm{Q}_r,\ \bm{f}^{(r)} = \bm{Q}_r^\trans\bm{f}.
\end{equation}
Once the reduced matrices $\bm{A}_i^{(r)}$, $\bm{B}_{ij}^{(r)}$, and the reduced vector $\bm{f}^{(r)}$ are computed, the linear system for each $ \sv\in\Gamma$ can be assembled with cost $O(r^2)$. It is important to note that $\vrbfemsoln$ represents  the coefficient vector associated with the basis~$\{q_1(\pv),\ldots,q_{r}(\pv)\}$. Therefore, the coefficient vector corresponding to the basis $\{v_1(\pv),\ldots, v_{\Nx}(\pv)\}$ can be obtained by multiplying  $\bm{Q}_r$ with $\vrbfemsoln$.

Typically, the reduced space is constructed using a greedy algorithm based on an indicator. In this study, the indicator is denoted by $\Delta_n(\sv)$ for the $n$-dimensional RBM space. The construction of $\bm{Q}_r$ using the greedy algorithm is provided by Algorithm~\ref{alg:RBM}.

\begin{algorithm}[h]
	\caption{Greedy algorithm for construction of reduced basis space}\label{alg:RBM}
	\begin{algorithmic}
		{\small
			\STATE {\bfseries Input:} a set of candidate parameters $\smp$, the RBM dimension $r$;
			\STATE  Randomly select $\sv^{(1)}\in \smp$ as the first sample, and set $n=1$;
			\STATE Compute $\bm{u}_{\sv^{(1)}}$ by solving \eqref{eq:linear_full} and set $\bm{Q}_1 = {\bm{u}_{\sv^{(1)}}}\big/\Vert {\bm{u}_{\sv^{(1)}}}\Vert_2$;
			\WHILE{ $n<r$}
			\FOR{each $\sv\in\smp$}
			\STATE Compute the indicator $\Delta_n(\sv)$;
			\ENDFOR
			\STATE $\sv^{(n+1)}=\arg\max\limits_{\sv\in\smp}\Delta_n(\sv)$;
			\STATE Compute the snapshots $\bm{u}_{\sv^{(n+1)}}$ by solving \eqref{eq:linear_full};
			\STATE Compute $\bm{q}$, the orthogonal complement of $\bm{u}_{\sv^{(n+1)}}$ with respect to $\bm{Q}_{n}$;
			\STATE Augment the reduced basis matrix $\bm{Q}_{n+1} = [\bm{Q}_n, \bm{q}]$ and set $n=n+1$.
			\ENDWHILE
		}
	\end{algorithmic}
\end{algorithm}

The most commonly used indicator in reduced basis methods is the residual based a posteriori error estimator. Detailed discussions are provided in \cite{Quarteroni2016Reduced}. However, in this work, where only the most important samples and their corresponding solutions are required, we employ a residual free indicator to generate the reduced basis space. We review the residual free indicator following the presentation in \cite{Chen2019resfree}.

Recall that the reduced solution, denoted as $\frbfemsoln(\pv,\sv)$, belongs to the span of the functions $\ffemsoln(\pv,\sv^{(1)}),\ldots,\ffemsoln(\pv,\sv^{(r)})$. Therefore, we can express it as follows:
\begin{equation}\notag
\frbfemsoln(\pv,\sv) = \sum_{i=1}^{r} u_h(\pv,\sv^{(i)}) l_i(\sv),
\end{equation}
where the expansion coefficients $l_i(\sv)$ are functions of $\sv$, and act as basis functions. These coefficients satisfy the following condition:
\begin{equation}\notag
l_i(\sv^{(j)}) = \delta_{ij},
\end{equation}
where $\delta_{ij}$ is the Kronecker delta. This means that the coefficient functions $l_j(\sv)$ are cardinal Lagrange interpolants associated with the space of functions defined by their span. To quantify the quality of the interpolation, we introduce the residual free indicator defined as:
\begin{equation}\notag
\Delta_n(\sv)= \sum_{i=1}^{n}|l_i(\sv)|.
\end{equation} 
The function above is essentially the norm of an interpolation operator and is referred to as the Lebesgue function in interpolation theory. Studies have demonstrated that the Lebesgue function captures the behavior of the residual based a posteriori error estimator, making it highly valuable in the selection of snapshots in the reduced basis method. Details of the residual free indicator are presented in \cite{Chen2019resfree}.

\subsection{The reduced basis  Galerkin method}
Suppose $Q_r$ is a reduced basis space generated by the RBM described in Section~\ref{sec:RBM}. We seek a reduced basis approximate solution $\frbsfemsoln(\pv,\sv)\in Q_r\otimes S_p$ that satisfies the following equation:
\begin{equation}\notag
\int_{\Gamma}\int_{D}\mypdf(\sv)\left[  a \nabla \frbsfemsoln \cdot \nabla w^*  - \kappa^2 \frbsfemsoln w^*\right]\rd \pv\rd \sv = \int_{\Gamma}\int_D \mypdf(\sv)fw^*\rd \pv\rd\sv,\ \forall  w \in Q_r\otimes S_p.
\end{equation}
Assuming that the reduced basis approximate solution $\frbsfemsoln(\pv,\sv)$ can be expanded as: 
\begin{equation}\notag
\frbsfemsoln(\pv,\sv)=\sum_{j=1}^{\Nxi}\sum_{s=1}^{r}u_{js}^{(r)}q_s(\pv)\Phi_j(\sv),
\end{equation}
where $u_{js}^{(r)}$ represents the coefficient of  $\frbsfemsoln(\pv,\sv)$ corresponding to $q_s(\pv)\Phi_j(\sv)$. To facilitate the presentation, we introduce the matrix $\mrbsfemsoln$ and the vector $\vrbsfemsoln$ as:
\begin{equation}\notag
\mrbsfemsoln(s,j)=u_{js}^{(r)}, \ \vrbsfemsoln = \mathtt{vec}(\mrbsfemsoln),\ \mrbsfemsoln \in \mathbb{R}^{r\times\Nxi}.
\end{equation}
By utilizing the relationship \eqref{eq:reducedbasis},  the coefficients of $\frbsfemsoln(\pv,\sv)$ can be obtained by solving the following reduced linear system:
\begin{equation}\label{eq:rbsgeq}
\left(\sum_{i=0}^{m}\bm{G}_{i0}\otimes\bm{A}_i^{(r)} -\sum_{i=0}^{m}\sum_{j=0}^{m}\bm{G}_{ij}\otimes\bm{B}_{ij}^{(r)}\right)\vrbsfemsoln =\bm{f}^{(r)},
\end{equation}
where 
\begin{equation}\notag
\bm{A}_{i}^{(r)} = \bm{Q}_r^\trans\bm{A}_i\bm{Q}_r, \ \bm{B}_{ij}^{(r)} = \bm{Q}_r^\trans\bm{B}_{ij}\bm{Q}_r,\ 
\bm{f}^{(r)} = \bm{Q}_r^\trans\bm{f}.
\end{equation}
Since $Q_r\subset V$, the reduced basis approximate solution $\frbsfemsoln(\pv,\sv)$ can also be expanded as:
\begin{equation}\notag
\frbsfemsoln(\pv,\sv)=\sum_{j=1}^{\Nxi}\sum_{s=1}^{\Nx}u_{js}v_s(\pv)\Phi_j(\sv),
\end{equation}
where the coefficients are given by
\begin{equation}\label{eq:rb2full}
u_{js} = \bm{U}(s,j), \ \bm{U} = \bm{Q}_r\mrbsfemsoln.  \ 
\end{equation}
It is worth noting that \eqref{eq:rbsgeq} has the same form as \eqref{eq:linear_sg}, but with a much smaller size. As a result, it can be solved much more efficiently.

Now, let us consider the residual, i.e., $\Vert \bm{b} - \bm{A}\vrbsfemsoln\Vert_2$. By \eqref{eq:linear_lhs}, \eqref{eq:mv2mm} and \eqref{eq:rb2full}, we have 
\begin{equation}\label{eq:res}
\Vert \bm{b} - \bm{A}\vrbsfemsoln\Vert_2  = \left\Vert \bm{b} - \mathtt{vec}\left(\sum_{i=0}^{m}(\bm{A}_i\bm{Q}_r)(\mrbsfemsoln\bm{G}_{i0}^\trans)+\sum_{i=0}^m\sum_{j=0}^m (\bm{B}_{ij}\bm{Q}_r)(\mrbsfemsoln\bm{G}_{ij}^\trans)\right)\right\Vert_2.  
\end{equation}
It is important to note that  $\bm{A}_{i}$, $\bm{B}_{ij}$, and $\bm{G}_{ij}$ are typically sparse matrices. Since $\bm{Q}_r\in R^{\Nx\times r}$, $\mrbsfemsoln\in R^{r\times \Nxi}$, the computational cost of evaluating \eqref{eq:res} can be estimated as:
\begin{equation}\notag
O\left(  (m+1)(m+2)[\Nx r+ \Nxi r+\Nxi\Nx r]\right) \approx O(m^2\Nx\Nxi r).
\end{equation}

Based on  the above descriptions, we propose a reduced basis stochastic Galerkin method, which is referred to as RBSGM in the following. Since the required number of reduced basis functions $r$ is unknown, the reduced basis space is constructed adaptively in RBSGM. The relative residual
\begin{equation}\label{eq:relres}
\mathtt{relres}=\Vert \bm{b}-\bm{A}\vrbsfemsoln\Vert_2/\Vert\bm{b}\Vert_2
\end{equation}
can be regarded as a function of the number of reduced basis functions $r$. To find an appropriate value of $r$ satisfying
\begin{equation}\label{eq:resleqtol}
\Vert \bm{b}-\bm{A}\vrbsfemsoln\Vert_2/\Vert\bm{b}\Vert_2 \leq \tol,
\end{equation}
where $\tol$ is the desired tolerance, we define a function $h(r)$ as the logarithm (base $10$) of the relative residual $\mathtt{relres}$, i.e.,
\begin{equation}\notag
h(r) = \lg(\mathtt{relres}), \ r = 1,2,\ldots,
\end{equation}
Obviously, \eqref{eq:resleqtol} is equivalent to  
\begin{equation}\label{eq:findr}
h(r) \leq \lg(\tol).
\end{equation}
To minimize  the number of times to compute $\mathtt{relres}$ using equation \eqref{eq:res}, it is desirable to find the value of $r$ satisfying \eqref{eq:findr} in as few steps as possible. In this work, we employ the secant method to find $r$ satisfying~\eqref{eq:findr}. We first compute $h(r_1)$ and $h(r_2)$, and the new  value of $r$ is determined by
\begin{equation}\label{eq:findr2}
r= r_1+\frac{r_2-r_1}{h(r_2)-h(r_1)}\left[\lg(\tol)-h(r_1)\right].
\end{equation}

The reduced basis stochastic Galerkin method is summarized in Algorithm~\ref{alg:RBGM}. It should be noted that, to ensure the robustness and efficiency of the algorithm, we increase the number of reduced basis functions by a multiple of $\mathtt{ns}$ during each adaptive iteration. In our approach, we refer to the addition of $\mathtt{ns}$ reduced basis functions as a stage, where $\mathtt{ns}$ represents the number of reduced basis functions added in each stage.
	
Suppose that the new value of $r$ is predicted using $r_1$, $r_2$, and \eqref{eq:findr2}. The number of stages $\mathtt{st}$ can be determined as follows:
\begin{equation}\notag
\mathtt{st}=\mathrm{floor}\left((r-r_2)/\mathtt{ns}\right)+1.
\end{equation}
Therefore, in the next iteration, $\mathtt{st}\cdot\mathtt{ns}$ reduced basis functions will be added to the reduced basis space. We use $k$ to denote the number of reduced basis functions that have already been added in the inner iteration. The inner iteration continues as long as $k< \mathtt{st}\cdot\mathtt{ns}$.

\begin{algorithm}[h]
	\caption{The reduced basis stochastic Galerkin method}\label{alg:RBGM}
	\begin{algorithmic}
		{\small
			\STATE {\bfseries Input:} The matrices $\bm{G}_{ij}$, $\bm{A}_{i}$, $\bm{B}_{ij}$ and  vectors $\bm{h}$, $\bm{f}$, a set of candidate parameters $\smp$, the RBM dimension $\mathtt{ns}$ in each stage, the maximum RBM dimension $\mathtt{nmax}$ in total, the tolerance of the iterative method $\tol$;
			\STATE Randomly select $\sv^{(1)}\in \smp$ as the first sample, set $n=1$, and set $\bm{Q}_1 = {\bm{u}_{\sv^{(1)}}}\big/\Vert {\bm{u}_{\sv^{(1)}}}\Vert_2$;
			\STATE Solve \eqref{eq:rbsgeq}, compute the relative residual $\mathtt{relres}$ via \eqref{eq:res}--\eqref{eq:relres}, and set $\mathtt{st}=1$, $r_1=1$;
			\WHILE{  $n<\mathtt{nmax}$ \& $\mathtt{relres}> \tol$}
			\STATE Set $k=0$;
			\WHILE{$n<\mathtt{nmax}$ \& $k<\mathtt{st}\cdot\mathtt{ns}$}			
			\FOR{each $\sv\in\smp$}
			\STATE Compute the indicator $\Delta_n(\sv)$;
			\ENDFOR
			\STATE $\sv^{(n+1)}=\arg\max\limits_{\sv\in\smp}\Delta_n(\sv)$;
			\STATE Compute the snapshots $\bm{u}_{\sv^{(n+1)}}$ by solving \eqref{eq:linear_full};
			\STATE Compute $\bm{q}$, the orthogonal complement of $\bm{u}_{\sv^{(n+1)}}$ with respect to $\bm{Q}_n$;
			\STATE Augment the reduced basis matrices 
			\begin{equation}\notag
			\begin{split}
			\bm{Q}_{n+1} &= [\bm{Q}_n, \bm{q}], \quad
			\bm{Q}_{n+1}^\trans \bm{f} = [\bm{Q}_n^\trans \bm{f}; \bm{q}^\trans \bm{f}],	\\
			\bm{A}_i\bm{Q}_{n+1} &= [\bm{A}_i\bm{Q}_n,\bm{A}_i\bm{q}],\quad  
			\bm{B}_{ij}\bm{Q}_{n+1}  = [\bm{B}_{ij}\bm{Q}_n,\bm{B}_{ij}\bm{q}], \\
			\bm{Q}_{n+1}^\trans \bm{A}_i\bm{Q}_{n+1} 
			&= \begin{bmatrix}
			\bm{Q}_n^\trans \bm{A}_i\bm{Q}_n&\bm{Q}_n^\trans \bm{A}_i\bm{q}\\
			\bm{q}^\trans \bm{A}_i\bm{Q}_n & \bm{q}^\trans\bm{A}_i\bm{q}
			\end{bmatrix}, \quad
			\bm{Q}_{n+1}^\trans \bm{B}_{ij}\bm{Q}_{n+1} 
			= \begin{bmatrix}
			\bm{Q}_n^\trans \bm{B}_{ij}\bm{Q}_n&\bm{Q}_n^\trans \bm{B}_{ij}\bm{q}\\
			\bm{q}^\trans \bm{B}_{ij}\bm{Q}_n & \bm{q}^\trans\bm{B}_{ij}\bm{q}
			\end{bmatrix};
			\end{split}
			\end{equation}
            \STATE Set $n=n+1$, $k=k+1$;
			\ENDWHILE
			\STATE Solve \eqref{eq:rbsgeq}, and compute the relative residual $\mathtt{relres}$ via \eqref{eq:res}--\eqref{eq:relres};
			\STATE Set $r_2=n$, predict $r$ via \eqref{eq:findr2}, set $r_1=r_2$, and  set $\mathtt{st}= \mathrm{floor}\left((r-r_2)/\mathtt{ns}\right)+1$.
			\ENDWHILE
			\RETURN $\bm{u} = \mathtt{vec}\left(\bm{Q}_n\bm{U}^{(n)}\right)$.
		}
	\end{algorithmic}
\end{algorithm}

\section{Numerical study}
In this section, we consider two problems: a diffusion problem and a Helmholtz problem. In both problems, the approximation obtained by the stochastic Galerkin method (SGM) with a gPC order of $p=6$ serves as the reference solution. For test problem 1, we use the standard finite element method to generate the matrices~$\bm{A}_i$ and $\bm{B}_{ij}$, where $i,j\in \{0,1,\ldots,m\}$. However, for test problem 2, we apply the codes associated with \cite{Liu2015Additive} for convenience. All the results presented in this study are obtained using MATLAB R2015b on a desktop computer with 2.90GHz Intel Core i7-10700 CPU.

\subsection{Test problem 1}
In this test problem, we consider the diffusion equation:
\begin{equation}\notag
\begin{cases}
-\nabla\cdot(a(\pv,\sv)\nabla u)=1, & \forall (\pv,\sv)\in D \times\Gamma\\
u = 0,   & \forall (\pv,\sv) \in \partial D \times \Gamma
\end{cases},
\end{equation} 
where $D = [-1,1]^2$ and
\begin{equation}\notag
a(\pv,\sv) = \mu + \sigma\sum_{i=1}^{m}\sqrt{\lambda_i}a_i(\pv)\svc_i
\end{equation} 
is a truncated KL expansion of random field  with 
mean function $\mu=0.2$, standard deviation $\sigma=0.1$ and 
covariance function 
\begin{equation}\notag
\mbox{cov}(\pv,\bm{y})=\sigma^2 \exp\left(-{|x_1-y_1|}-{|x_2-y_2|}
\right).
\end{equation}
The random variables $\svc_i$ are
chosen to be identically independent distributed uniform random variables on $[-1, 1].$

In this test problem, the gPC order is set to $5$, the number of reduced basis functions in each stage is set to $15$, and the number of candidate parameters is set to $500$. In Algorithm~\ref{alg:RBGM}, we need to solve the reduced linear system~\eqref{eq:rbsgeq} during each adaptive iteration using an iterative method. To ensure accuracy, the tolerance for the reduced linear system is typically chosen to be smaller than $\tol$. In this test problem, we set $\tol$ to $10^{-4}$ or $10^{-5}$, and the tolerance for the reduced linear system to $10^{-7}$. Both the full linear systems~\eqref{eq:linear_sg} and the reduced linear systems~\eqref{eq:rbsgeq} are solved using the preconditioned conjugate gradients method (PCG) with a mean based preconditioner.

\begin{table}[htbp]
	\caption{CPU time for different $m$ and $\Nx$ with $\mathtt{ns}=15$ and $p=5$.}\label{tab:test1a}
	\begin{center}
	\newcolumntype{C}{>{\centering\arraybackslash}X}%
	\begin{tabularx}{0.95\linewidth}{c|CCCC}
		\hline
		$\tol$&	$\Nx$ & $m = 5$  & $m=7$ & $m=10$ \\
		\hline
		\multirow{3}*{$10^{-4}$} 	& $33^2$&  4.67 (0.42) &  12.57 (1.70) &  40.89 (9.05)  \\  
		& $65^2$&  5.97 (2.41) &  15.39 (10.74) &  52.79 (56.38) \\  
		& $129^2$& 11.63 (19.14) &  31.12 (93.41) &  76.18 (417.95) \\ 
		\hline
		\multirow{3}*{$10^{-5}$}	& $33^2$& 13.11 (0.58) &  28.55 (2.52) &  75.73 (12.40) \\  
		& $65^2$& 15.72 (3.30) &  37.31 (13.61) &  98.76 (77.58)\\  
		& $129^2$& 25.49 (26.82) &  63.03 (120.27) &  162.88 (580.99) \\
		\hline
	\end{tabularx}
\end{center}
\end{table}
In Table~\ref{tab:test1a}, we present the CPU time of RBSGM for different values of $m$ and $\Nx$. Additionally, we provide the CPU time of SGM for the same tolerance, indicated within the brackets. The table clearly demonstrates that RBSGM outperforms SGM in terms of efficiency when the dimension of the stochastic space is high or when the grid in the physical space is sufficiently fine.

\begin{table}[htbp]
	\caption{CPU time for generating the reduced basis functions and their percentages of the total CPU time.}\label{tab:test1b}
\begin{center}
	\newcolumntype{C}{>{\centering\arraybackslash}X}%
	\begin{tabularx}{0.95\linewidth}{c|CCCC}
		\hline
		$\tol$ & $\Nx$  & $m = 5$  & $m=7$ & $m=10$ \\
		\hline
		\multirow{3}*{$10^{-4}$} 	& $33^2$&  4.53 (97) &  11.74 (93) &  30.07 (74) \\  
		& $65^2$&  5.83 (98) &  14.25 (93) &  39.23 (74) \\  
		& $129^2$& 11.19 (96) &  28.70 (92) &  60.51 (79) \\ 
		\hline
		\multirow{3}*{$10^{-5}$} 	& $33^2$& 12.82 (98) &  27.18 (95) &  57.79 (76) \\  
		& $65^2$& 15.46 (98) &  35.10 (94) &  76.20 (77) \\  
		&$129^2$& 24.94 (98) &  59.07 (94) &  129.78 (80) \\  
		\hline
	\end{tabularx}
\end{center}	
\end{table}

In Table~\ref{tab:test1b}, we provide the CPU time required for generating the reduced basis functions, along with their percentages of the total CPU time indicated in the brackets. It is evident that in RBSGM, the computation for generating the reduced basis functions consumes the majority of the CPU time, making it the main bottleneck of the proposed method in terms of efficiency. However, the development of an efficient method to generate the reduced basis functions is beyond the scope of this work.

\begin{figure}[htbp]
	\begin{center}
		\subfloat[$m=5,N_h=33^2$]{
			\includegraphics[width=0.32\linewidth]{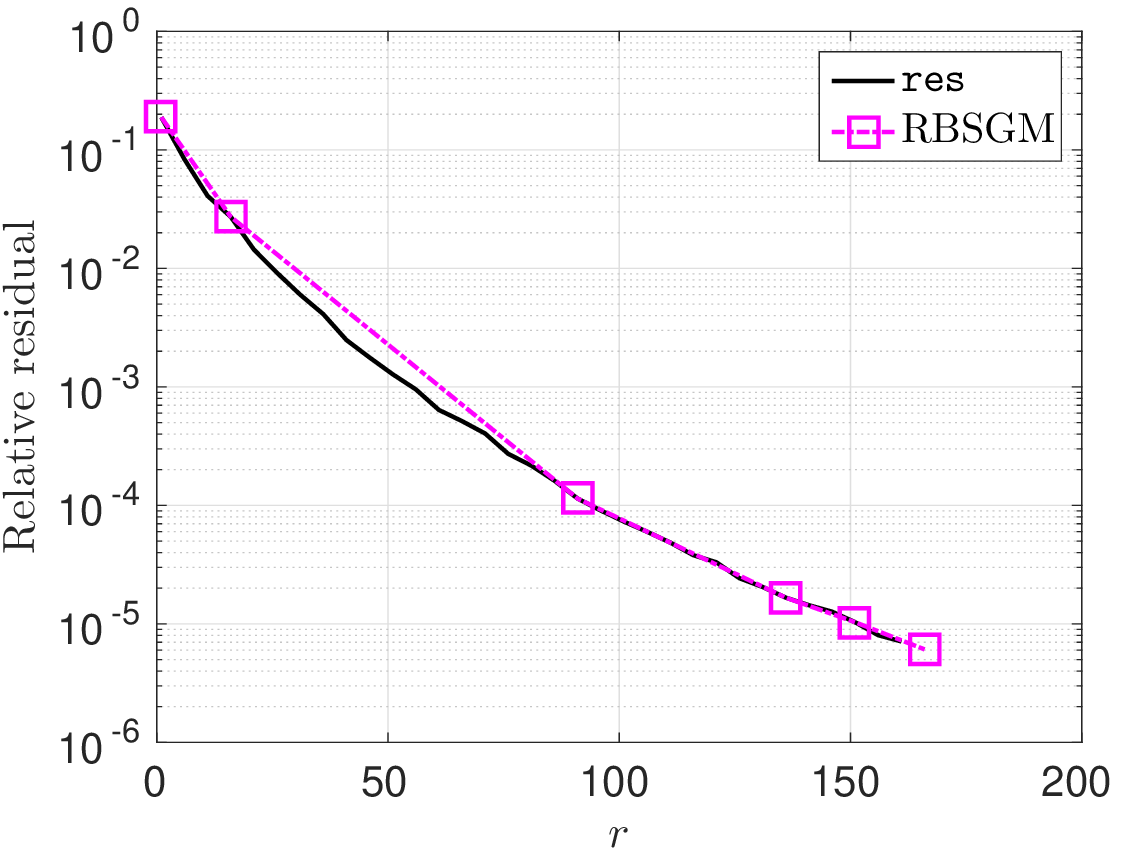}
		}
		\subfloat[$m=5,N_h=65^2$]{
			\includegraphics[width=0.32\linewidth]{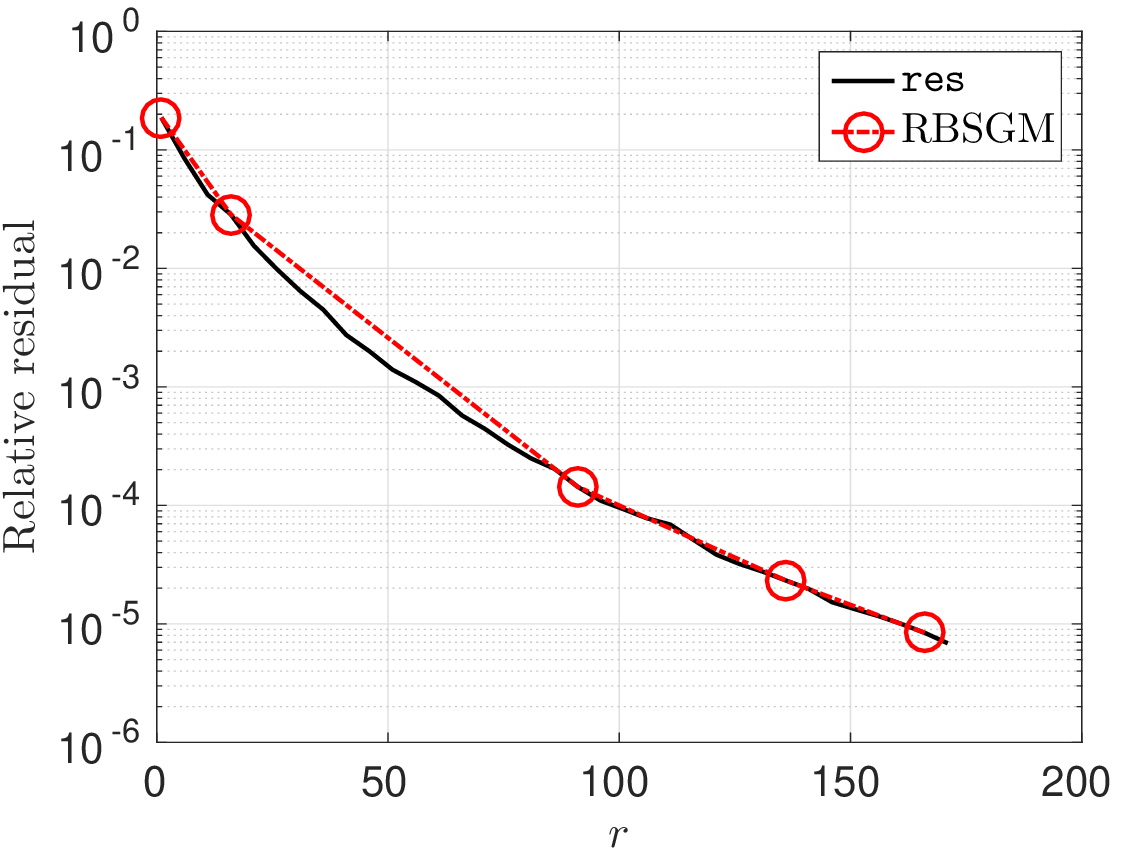}
		}
		\subfloat[$m=5,N_h=129^2$]{
			\includegraphics[width=0.32\linewidth]{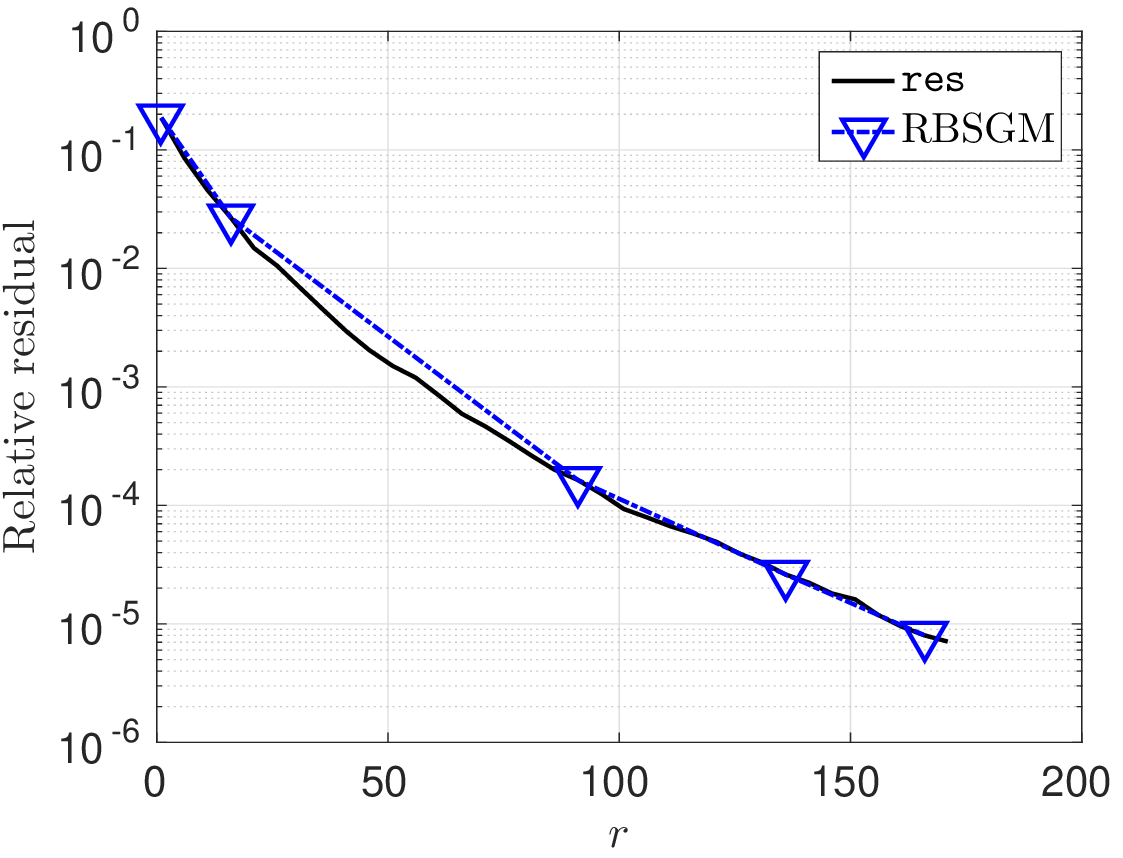}
		}\\
		\subfloat[$m=7,N_h=33^2$]{
			\includegraphics[width=0.32\linewidth]{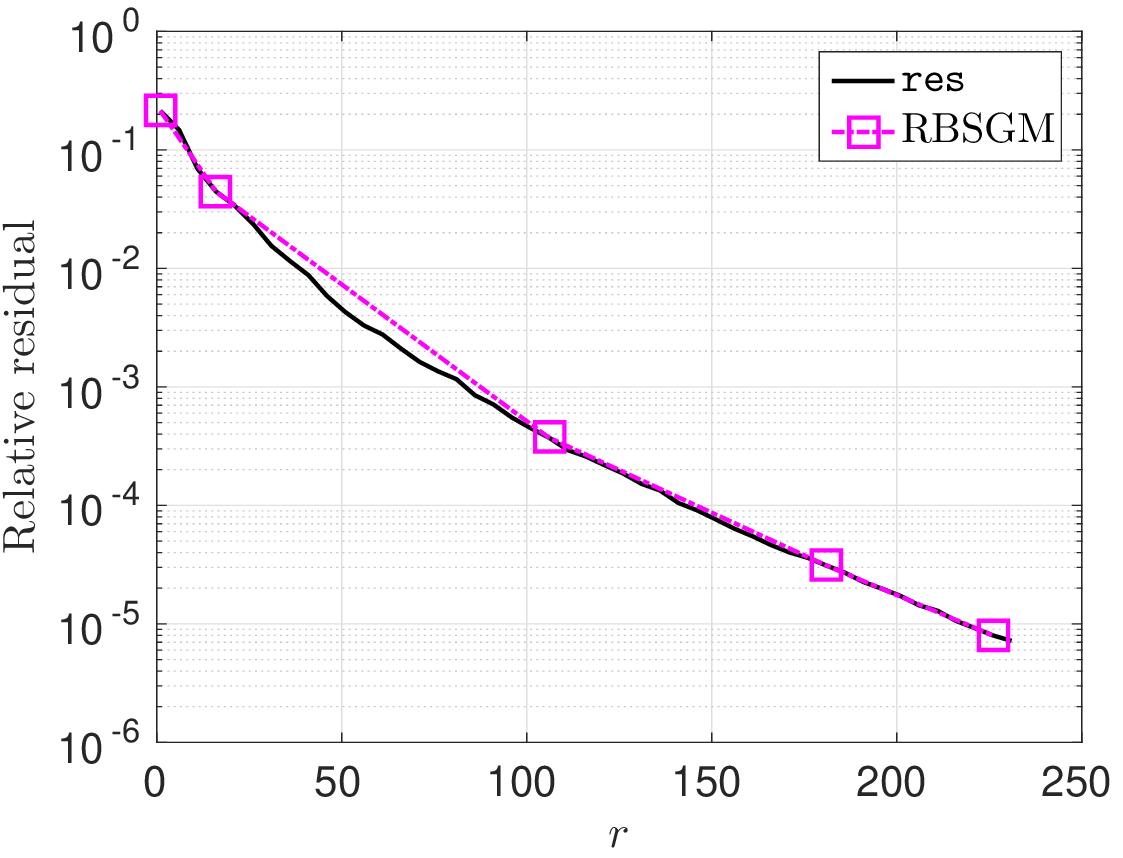}
		}
		\subfloat[$m=7,N_h=65^2$]{
			\includegraphics[width=0.32\linewidth]{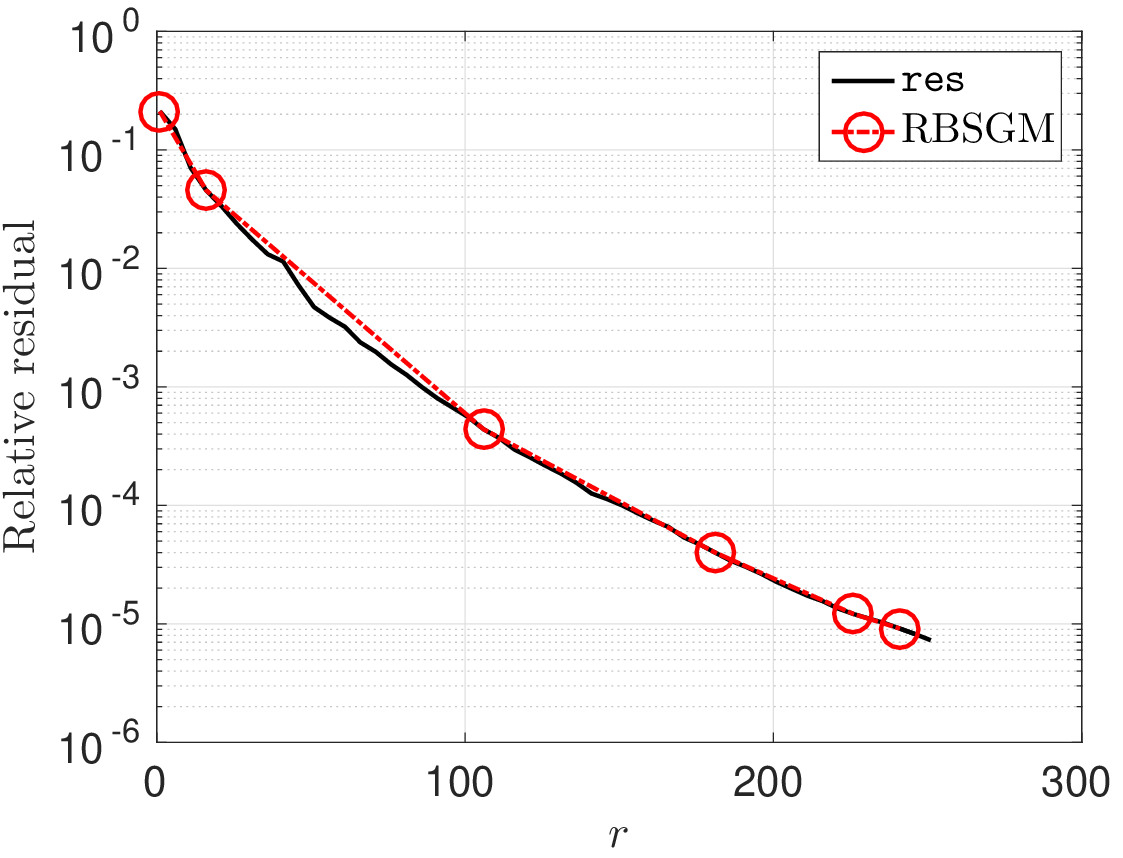}
		}
		\subfloat[$m=7,N_h=129^2$]{
			\includegraphics[width=0.32\linewidth]{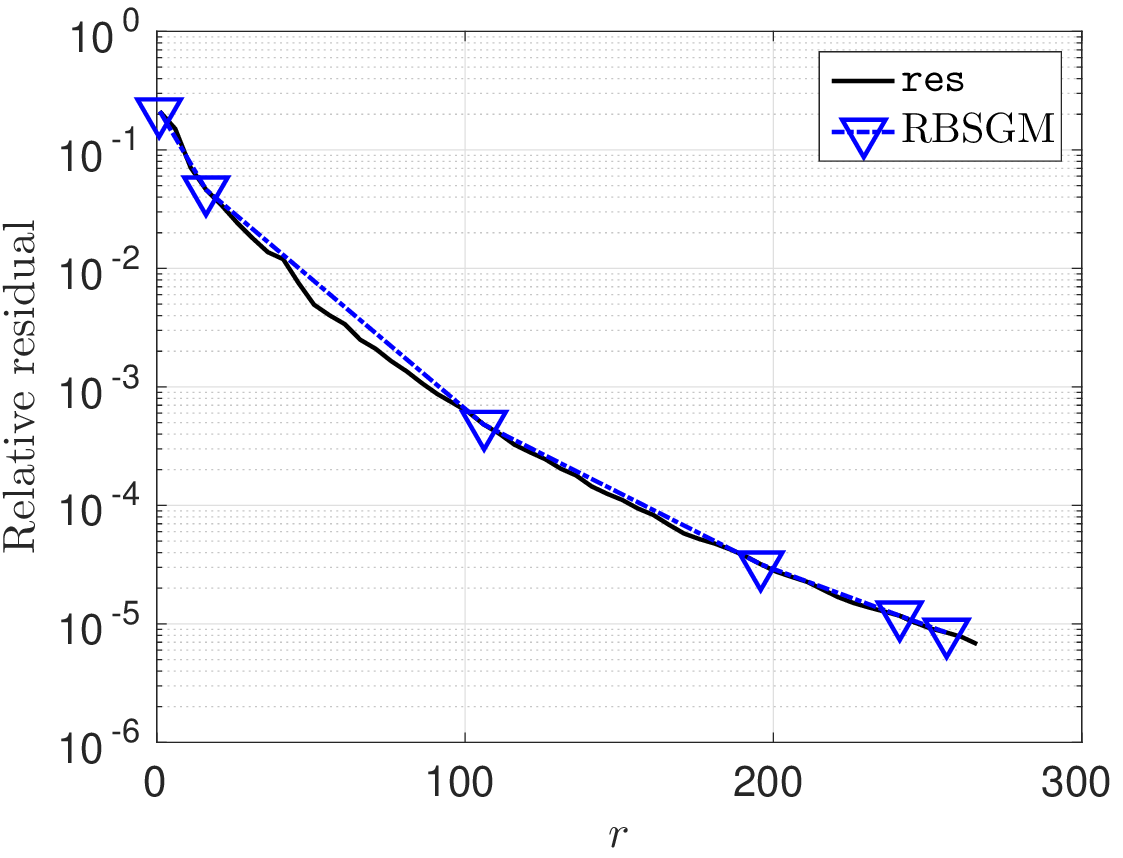}
		}\\
		\subfloat[$m=10,N_h=33^2$]{\includegraphics[width=0.32\linewidth]{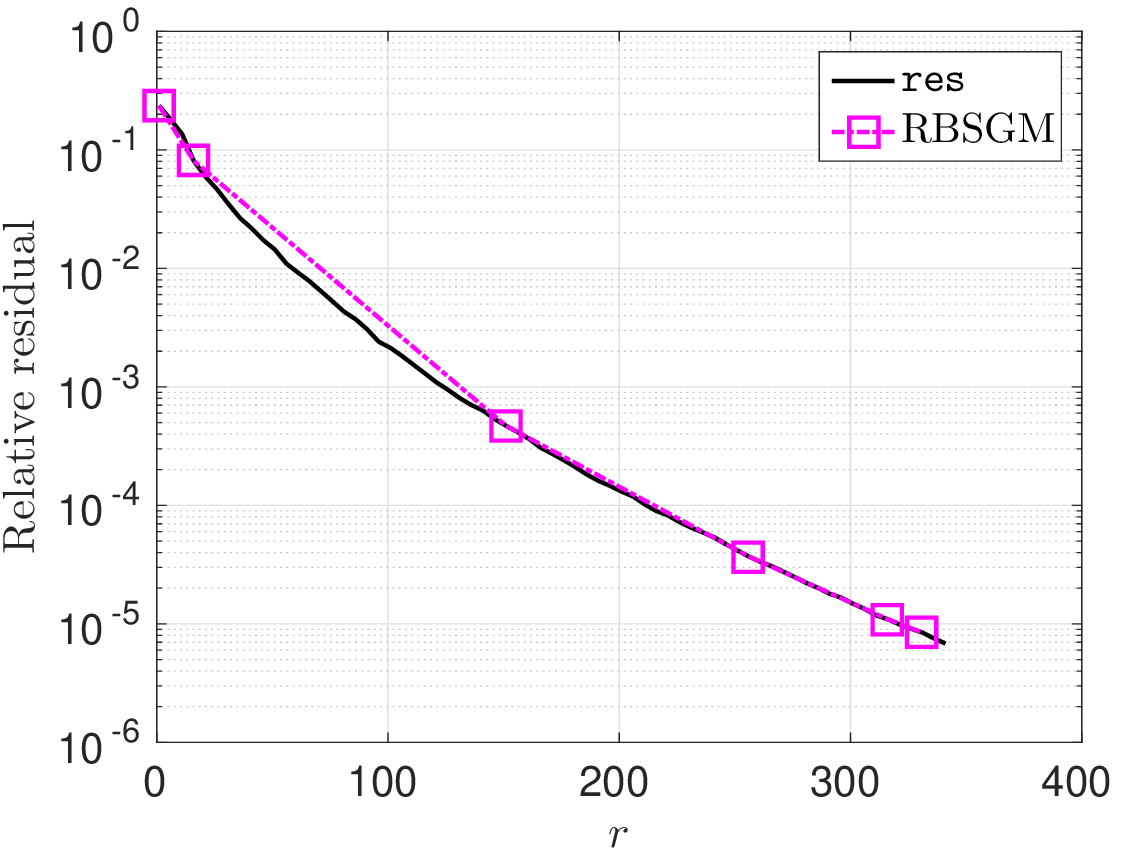}
		}
		\subfloat[$m=10,N_h=65^2$]{
			\includegraphics[width=0.32\linewidth]{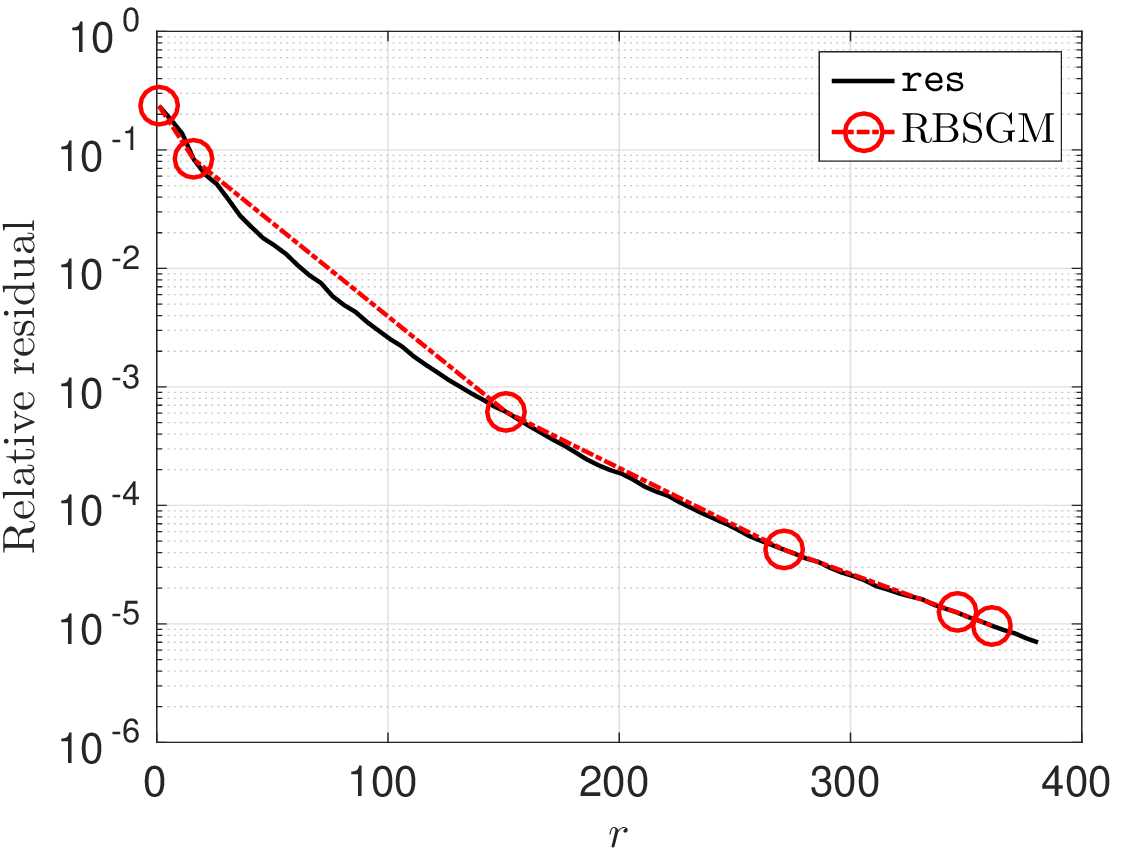}
		}
		\subfloat[$m=10,N_h=129^2$]{
			\includegraphics[width=0.32\linewidth]{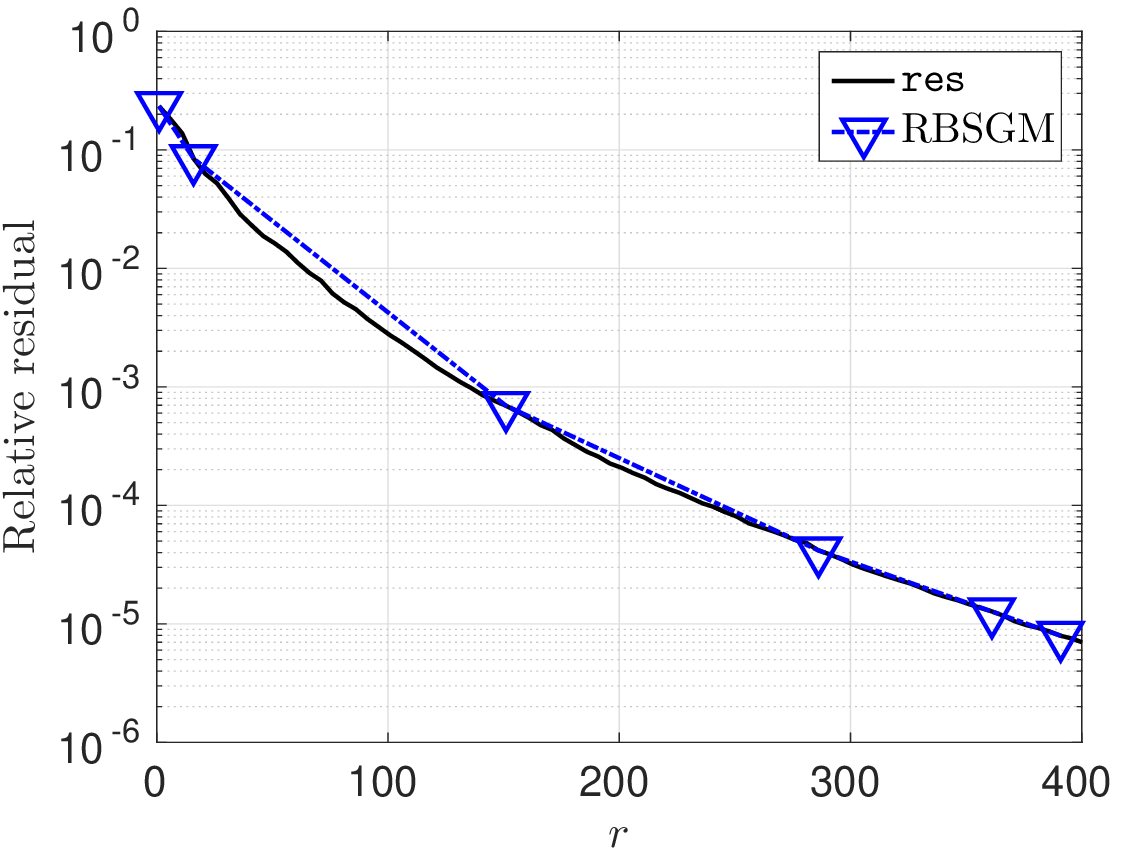}
		}\\
	\end{center}
	\caption{Relative residual with respect to the number of reduced basis functions $r$.}\label{fig:test1a}
\end{figure}

In Figure~\ref{fig:test1a}, we present the relationship between the number of reduced basis functions and the corresponding relative residual. Additionally, we provide the relative residual computed during the adaptive procedure of the RBSGM. Both the full linear system~\eqref{eq:linear_sg} and the reduced linear system~\eqref{eq:rbsgeq} are solved using the PCG with a mean based preconditioner. The figure clearly demonstrates that the number of reduced basis functions $r$ is significantly smaller than the DOF in the finite element method. Furthermore, we observe that $r$ is primarily influenced by $\tol$ but remains almost independent of~$\Nx$. In the RBSGM, the relative residual only needs to be computed $5$ or $6$ times when $\tol = 10^{-5}$. This highlights the computational efficiency of the reduced basis stochastic Galerkin method, as it achieves comparable accuracy with significantly smaller computational costs.

To access the accuracy of the reduced basis stochastic Galerkin approach, we consider the relative errors of the mean and variance functions, which are defined through
\begin{equation}\label{eq:error_define}
\meanerr = \frac{\Vert\pE{\frbsfemsoln}-\pE{\frefsoln}\Vert_{L^2}}{\Vert\pE{\frefsoln}\Vert_{L^2}},\ \varerr = \frac{\Vert\pV{\frbsfemsoln}-\pV{\frefsoln}\Vert_{L^2}}{
	\Vert\pV{\frefsoln}\Vert_{L^2}},
\end{equation}
where $\frefsoln$ is the reference solution with the mesh size $h$ in the physical domain. It is important to note that the reference solution varies with different mesh sizes. Similarly, we access the errors of the SGM using~\eqref{eq:error_define}, with $\frbsfemsoln$ replaced by $\fsfemsoln$.  In this test problem, we use the approximation $\fsfemsoln(\pv,\sv)$ obtained by the stochastic Galerkin method with gPC order $p=6$ as the reference solution. The reference solution is computed using the PCG with a mean based preconditioner and a tolerance of~$10^{-7}$.

\begin{table}[htbp]
	\caption{Relative errors of the mean and variance functions.}\label{tab:test1c}
\begin{center}
	\newcolumntype{C}{>{\centering\arraybackslash}X}%
	\begin{tabularx}{0.95\linewidth}{c|c|cCCC}
		\hline
		$\tol$ &  &$\Nx$ & $m = 5$  & $m=7$ & $m=10$\\
		\hline
		\multirow{6}*{$10^{-4}$}&\multirow{3}*{$\meanerr$}  
		    &$33^2$& 4.44e-07 (7.76e-07)  & 1.50e-06 (7.34e-06)  & 3.97e-06 (2.34e-06)  \\  
		 &  &$65^2$& 4.41e-07 (7.74e-07)  & 1.51e-06 (1.51e-06)  & 3.98e-06 (2.34e-06)  \\  
		 &  &$129^2$& 4.36e-07 (7.73e-07)  & 1.51e-06 (1.51e-06)  & 3.99e-06 (2.34e-06)  \\  
		 \cline{2-6}
		 & \multirow{3}*{$\varerr$} 
		    &$33^2$& 4.40e-05 (4.85e-05)  & 1.38e-04 (1.38e-04)  & 3.10e-04 (3.13e-04)  \\  
		 &  &$65^2$& 4.40e-05 (4.84e-05)  & 1.38e-04 (1.46e-04)  & 3.10e-04 (3.13e-04)  \\  
		 &  &$129^2$& 4.39e-05 (4.83e-05)  & 1.38e-04 (1.46e-04)  & 3.10e-04 (3.13e-04)  \\  
		 \hline
		\multirow{6}*{$10^{-5}$}&\multirow{3}*{$\meanerr$}  
		   & $33$& 4.18e-07 (4.26e-07)  & 1.50e-06 (1.49e-06)  & 3.97e-06 (3.97e-06)  \\  
		&  & $65$& 4.19e-07 (4.27e-07)  & 1.50e-06 (1.50e-06)  & 3.98e-06 (3.99e-06)  \\  
		&  & $129$& 4.20e-07 (4.28e-07)  & 1.51e-06 (1.50e-06)  & 3.99e-06 (3.99e-06)  \\  
		\cline{2-6}
		&\multirow{3}*{$\varerr$}  
		   &$33^2$& 4.39e-05 (4.72e-05)  & 1.38e-04 (1.37e-04)  & 3.10e-04 (3.10e-04)  \\  
		&  &$65^2$& 4.39e-05 (4.72e-05)  & 1.38e-04 (1.37e-04)  & 3.10e-04 (3.10e-04)  \\  
		&  &$129^2$& 4.39e-05 (4.72e-05)  & 1.38e-04 (1.37e-04)  & 3.10e-04 (3.10e-04)  \\  
		\hline
	\end{tabularx}
\end{center}	
\end{table}

Table~\ref{tab:test1c} presents the relative errors of the mean and variance functions for RBSGM. In addition, we provide the relative errors of the SGM in parentheses for comparison. The table reveals that the relative errors of the mean and variance functions in RBSGM are comparable to those in the SGM.  It is worth noting that, in the adaptive procedure, we solve the reduced basis linear system~\eqref{eq:rbsgeq} using the PCG with a tolerance of  $10^{-7}$. This tolerance value is smaller than $\tol$ used in this test problem. Meanwhile, since we increase the number of reduced basis functions by a multiple of $\mathtt{ns}=15$ during each adaptive iteration in Algorithm~\ref{alg:RBGM}, the number of reduced basis functions determined by the adaptive procedure may be slightly larger than that required. As a result, RBSGM outperforms SGM in terms of accuracy for many cases in this test problem.

\begin{table}[htbp]
	\caption{Overall relative errors of the mean and variance functions.}\label{tab:test1d}
		\begin{center}
			\newcolumntype{C}{>{\centering\arraybackslash}X}%
			\begin{tabularx}{0.95\linewidth}{c|c|cCCC}
				\hline
				$\tol$ &  &$\Nx$ & $m = 5$  & $m=7$ & $m=10$\\
				\hline
				\multirow{6}*{$10^{-4}$}&\multirow{3}*{$\meanerr^{\mathrm{o}}$}  
				   &$33^2$& 1.46e-03 (1.46e-03)  & 1.47e-03 (1.47e-03)  & 1.48e-03 (1.48e-03)  \\  
				&  &$65^2$& 3.70e-04 (3.70e-04)  & 3.72e-04 (3.72e-04)  & 3.76e-04 (3.75e-04)  \\  
				&  &$129^2$& 9.43e-05 (9.43e-05)  & 9.52e-05 (9.52e-05)  & 9.71e-05 (9.60e-05)  \\  
				\cline{2-6}
				&\multirow{3}*{$\varerr^{\mathrm{o}}$}   
				   &$33^2$& 3.04e-03 (3.04e-03)  & 3.16e-03 (3.16e-03)  & 3.23e-03 (3.23e-03)  \\  
				&  &$65^2$& 7.55e-04 (7.57e-04)  & 7.78e-04 (7.77e-04)  & 8.06e-04 (8.06e-04)  \\  
				&  &$129^2$& 1.84e-04 (1.87e-04)  & 2.16e-04 (2.20e-04)  & 3.33e-04 (3.36e-04)  \\  
				\hline
				\multirow{6}*{$10^{-5}$}&\multirow{3}*{$\meanerr^{\mathrm{o}}$}  
				   &$33^2$& 1.46e-03 (1.46e-03)  & 1.47e-03 (1.47e-03)  & 1.48e-03 (1.48e-03)  \\  
				&  &$65^2$& 3.70e-04 (3.70e-04)  & 3.72e-04 (3.72e-04)  & 3.76e-04 (3.76e-04)  \\  
				&  &$129^2$& 9.43e-05 (9.43e-05)  & 9.52e-05 (9.52e-05)  & 9.71e-05 (9.71e-05)  \\  
				\cline{2-6}
				&\multirow{3}*{$\varerr^{\mathrm{o}}$}  
				   &$33^2$& 3.04e-03 (3.04e-03)  & 3.16e-03 (3.16e-03)  & 3.23e-03 (3.23e-03)  \\  
				&  &$65^2$& 7.55e-04 (7.54e-04)  & 7.78e-04 (7.78e-04)  & 8.06e-04 (8.06e-04)  \\  
				&  &$129^2$& 1.84e-04 (1.84e-04)  & 2.16e-04 (2.16e-04)  & 3.33e-04 (3.33e-04)  \\ 
				\hline
			\end{tabularx}
		\end{center}	
\end{table}

In Table~\ref{tab:test1c}, our main focus is to compare the errors between the reduced solution (solution of RBSGM) and the full solution (solution of SGM). The relative errors, as defined in~\eqref{eq:error_define}, primarily reflect the accuracy of the stochastic approximation. To further examine the overall errors, we define the overall relative errors of the mean and variance functions as follows:
\begin{equation}\label{eq:overallerror_define}
\meanerr^{\mathrm{o}} = \frac{\Vert\pE{\frbsfemsoln}-\pE{\ofrefsoln}\Vert_{L^2}}{\Vert\pE{\ofrefsoln}\Vert_{L^2}},\ \varerr^{\mathrm{o}} = \frac{\Vert\pV{\frbsfemsoln}-\pV{\ofrefsoln}\Vert_{L^2}}{
	\Vert\pV{\ofrefsoln}\Vert_{L^2}},
\end{equation}
where $\ofrefsoln$ represents the reference solution obtained by the SGM with a gPC order of $p=6$, a mesh size $\Nx = 257^2$ in the physical domain, and a tolerance of $10^{-7}$ in PCG. Similarly, we assess the overall errors of the SGM using~\eqref{eq:overallerror_define}, with $\frbsfemsoln$ replaced by $\fsfemsoln$.

Table~\ref{tab:test1d} presents the overall relative errors of the mean and variance functions for RBSGM. In addition, we provide the relative errors of the SGM in parentheses for comparison. The table reveals that the overall relative errors of the mean and variance functions in RBSGM are comparable to those in the SGM. It is evident that the overall relative errors decrease as the resolution of the physical space increases. Furthermore, it is worth noting that the overall relative errors for $\tol=10^{-4}$ coincide with the errors for $\tol=10^{-5}$ with an accuracy of at least $6$ digits. This observation suggests that $\tol=10^{-4}$ can provide a satisfactory level of accuracy, while decreasing the tolerance for PCG further would not lead to a significant increase in the overall accuracy.

\subsection{Test problem 2}
In this test problem, we consider a stochastic Helmholtz equation
\begin{equation}\notag
-\nabla^2u - \kappa^2(\pv,\sv)u=f(\pv), \ \forall(\pv,\sv)\in D\times\Gamma,
\end{equation}
with the Sommerfeld radiation boundary condition. Let the domain of interest be $D = [0,1]^2$ and 
\begin{equation}\notag
\kappa(\pv,\sv) = \mu + \sigma\sum_{i=1}^{m}\sqrt{\lambda_i}\kappa_i(\pv)\svc_i
\end{equation}
is a truncated KL expansion of random field  with 
mean function $\mu=4\cdot(2\pi)$, standard deviation $\sigma=0.1\mu$, and 
covariance function 
\begin{equation}\notag
\mbox{cov}(\pv,\bm{y})=\sigma^2 \exp\left(-{|x_1-y_1|/4}-{|x_2-y_2|/4}
\right).
\end{equation}
The random variables $\svc_i$ are
chosen to be identically independent distributed uniform random variables on $[-1, 1].$ The  Gaussian point source at the center of the domain is served as the source term, i.e.,
\begin{equation}\notag
f(\pv) = - \mbox{e}^{-(8\cdot4)^2((x_1-0.5)^2+(x_2-0.5)^2)}.
\end{equation}

In this example, we use the perfectly matched layers (PML) to simulate the Sommerfeld condition~\cite{Berenger94}, and apply the codes associated with  \cite{Liu2015Additive} to generate the matrices $\bm{A}_i$, $\bm{B}_{ij}$ and the snapshots in Algorithm~\ref{alg:RBM}. Furthermore, the gPC order is set to $5$, the number of reduced basis functions in each stage is set to $10$, and the number of candidate parameters is set to $400$. In Algorithm~\ref{alg:RBGM}, we need to solve the reduced linear system~\eqref{eq:rbsgeq} during each adaptive iteration using an iterative method. To ensure accuracy, the tolerance for the reduced linear system is typically chosen to be smaller than $\tol$. In this test problem, we set $\tol$ to $10^{-4}$ or $10^{-5}$, while the tolerance for the reduced linear system is set to $10^{-7}$. Both the full linear systems~\eqref{eq:linear_sg} and the reduced linear systems~\eqref{eq:rbsgeq} are solved using the bi-conjugate gradient stabilized method (Bi-CGSTAB) with a mean based preconditioner. Additionally, we use the approximation $\fsfemsoln(\pv,\sv)$ obtained by the stochastic Galerkin method with gPC order $p=6$ as the reference solution. The reference solution is computed using the Bi-CGSTAB with a mean based preconditioner and a tolerance of~$10^{-7}$. It is important to note that the reference solution varies with different mesh sizes.

\begin{table}[htbp]
	\caption{CPU time for different $m$ and $\Nx$ with $\mathtt{ns}=10$ and $p=5$.}\label{tab:test2a}
\begin{center}
	\newcolumntype{C}{>{\centering\arraybackslash}X}%
	\begin{tabularx}{0.95\linewidth}{c|CCCC}
		\hline
		$\tol$ & $\Nx$ & $m = 5$  & $m=7$ & $m=10$ \\
		\hline
		\multirow{3}*{$10^{-4}$}& $33^2$&  4.94 (0.94) &  14.84 (4.43) &  56.48 (24.51) \\  
		& $65^2$&  6.36 (6.01) &  39.34 (25.88) &  158.76 (189.03) \\  
		& $129^2$ & 17.97 (40.99) &  69.45 (213.13) &  417.23 (1123.50) \\ 
		\hline 
		\multirow{3}*{$10^{-5}$}	& $33^2$& 11.24 (1.11) &  43.15 (5.25) &  150.57 (29.17) \\  
		& $65^2$& 17.56 ( 6.96) &  87.36 (29.82) &  337.68 (246.95)\\  
		& $129^2$& 38.19 (55.98) &  162.84 (281.06) &  695.23 (1468.82) \\   
		\hline
	\end{tabularx}
\end{center}
\end{table}
Table~\ref{tab:test2a} presents the CPU times of RBSGM for different values of $m$ and $\Nx$, along with the CPU times of SGM indicated in parentheses for the same tolerance. The results clearly demonstrate that RBSGM outperforms SGM in terms of efficiency when the dimension of the stochastic space is large or when the grid in the physical space is sufficiently fine.

\begin{table}[htbp]
	\caption{CPU time for generating the reduced basis functions and their percentages of the total CPU time.}\label{tab:test2b}
\begin{center}
	\newcolumntype{C}{>{\centering\arraybackslash}X}%
	\begin{tabularx}{0.95\linewidth}{c|CCCC}
		\hline
		$\tol$ & $\Nx$ & $m = 5$  & $m=7$ & $m=10$ \\
		\hline
		\multirow{3}*{$10^{-4}$}		&$33^2$&  4.49 (91) &  11.75 (79) &  28.97 (51) \\  
		&$65^2$&  5.47 (86) &  26.00 (66) &  65.58 (41) \\  
		&$129^2$& 12.89 (72) &  37.86 (55) &  122.01 (29) \\  
		\hline
		\multirow{3}*{$10^{-5}$}		&$33^2$& 10.51 (94) &  36.82 (85) &  84.77 (56)\\  
		&$65^2$& 16.00 (91) &  67.04 (77) &  158.35 (47) \\  
		&$129^2$& 30.91 (81) &  115.52 (71) &  311.68 (45) \\   
		\hline
	\end{tabularx}
\end{center}
\end{table}

Table~\ref{tab:test2b} presents the CPU time for generating the reduced basis functions, with the percentages of the total CPU time provided in the brackets. In this problem, the number of non-zero matrices in~\eqref{eq:linear_lhs} is $(m+1)^2+1$, which is significantly larger than that of test problem 1. As a result, the percentages of the computational cost associated with generating the reduced basis functions decrease, but it still remains a substantial part of the overall computational cost in the proposed method.

\begin{figure}[!htbp]
	\begin{center}
		\subfloat[$m=5,N_h=33^2$]{
			\includegraphics[width=0.32\linewidth]{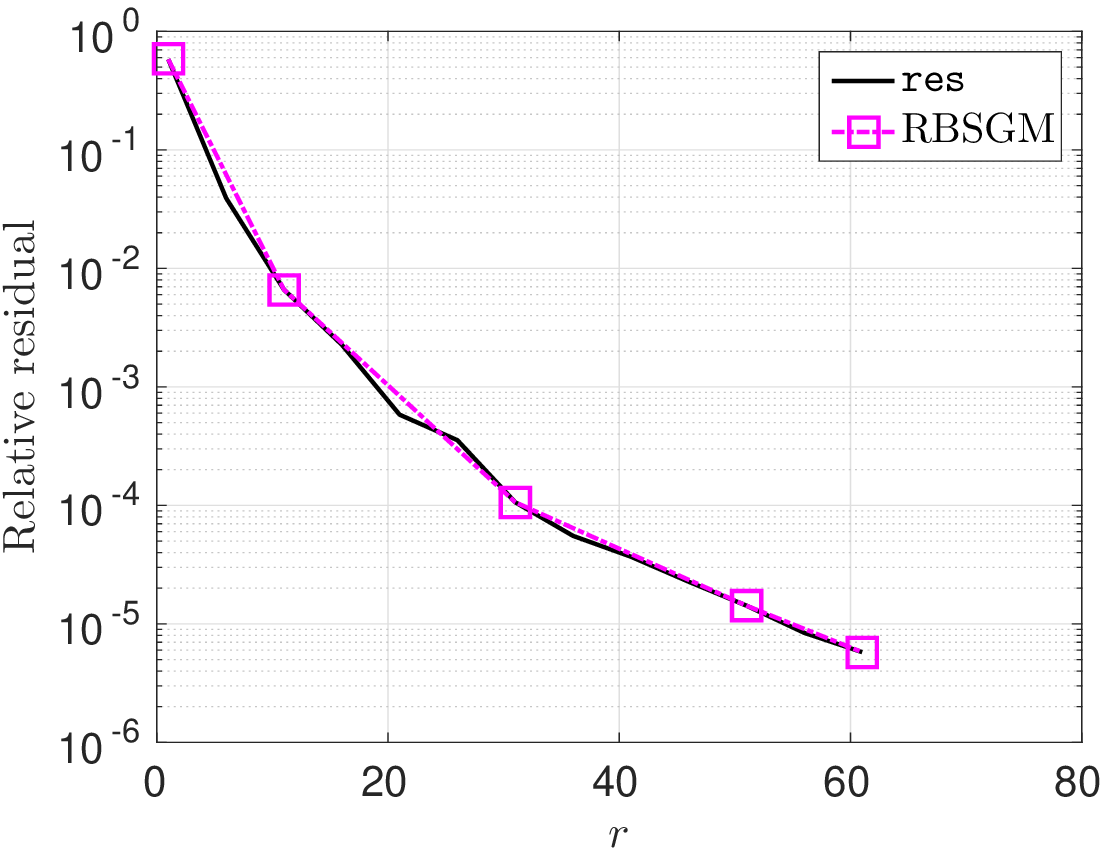}
		}
		\subfloat[$m=5,N_h=65^2$]{
			\includegraphics[width=0.32\linewidth]{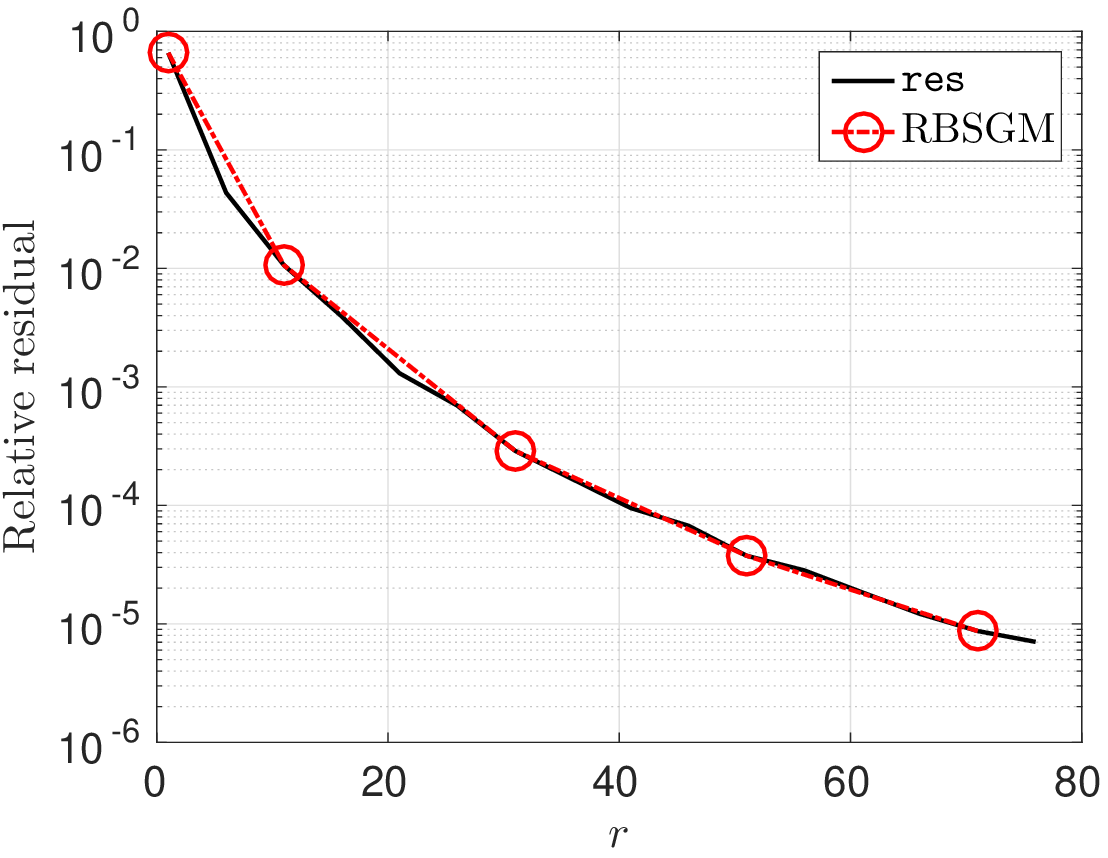}
		}
		\subfloat[$m=5,N_h=129^2$]{
			\includegraphics[width=0.32\linewidth]{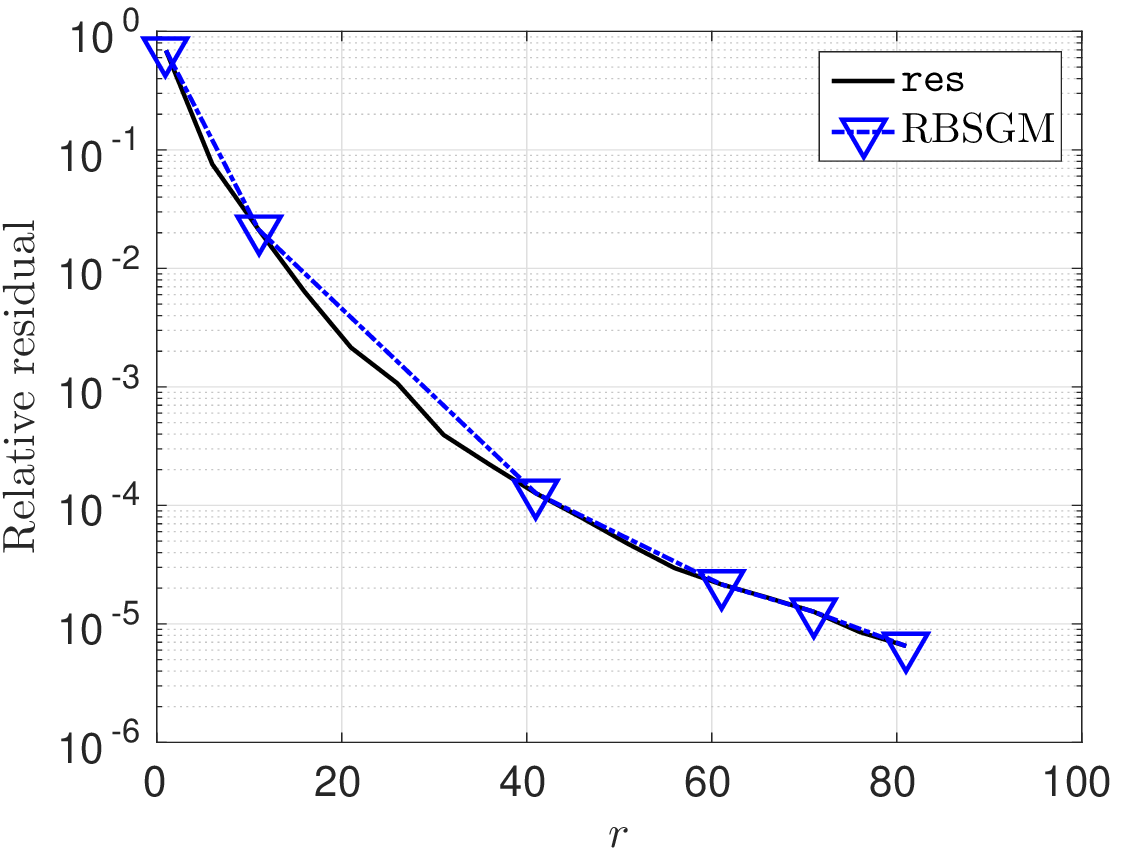}
		}\\
		\subfloat[$m=7,N_h=33^2$]{
			\includegraphics[width=0.32\linewidth]{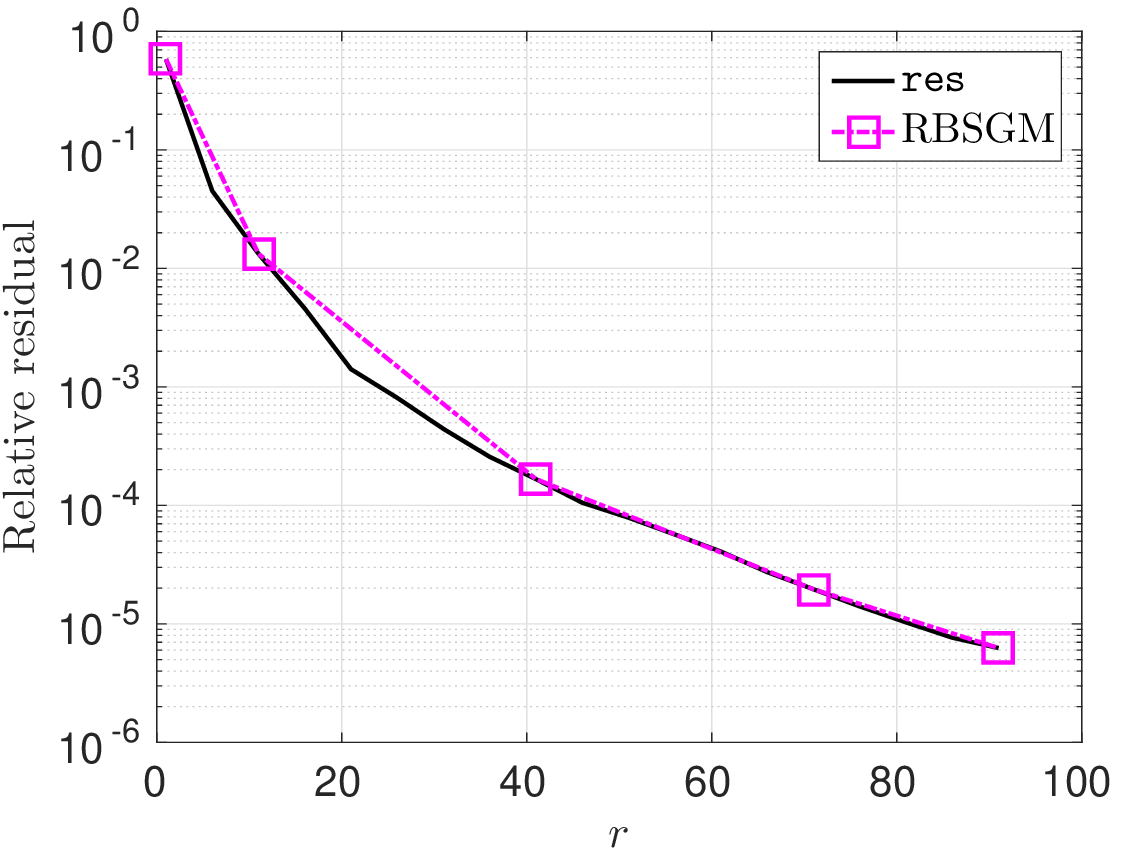}
		}
		\subfloat[$m=7,N_h=65^2$]{
			\includegraphics[width=0.32\linewidth]{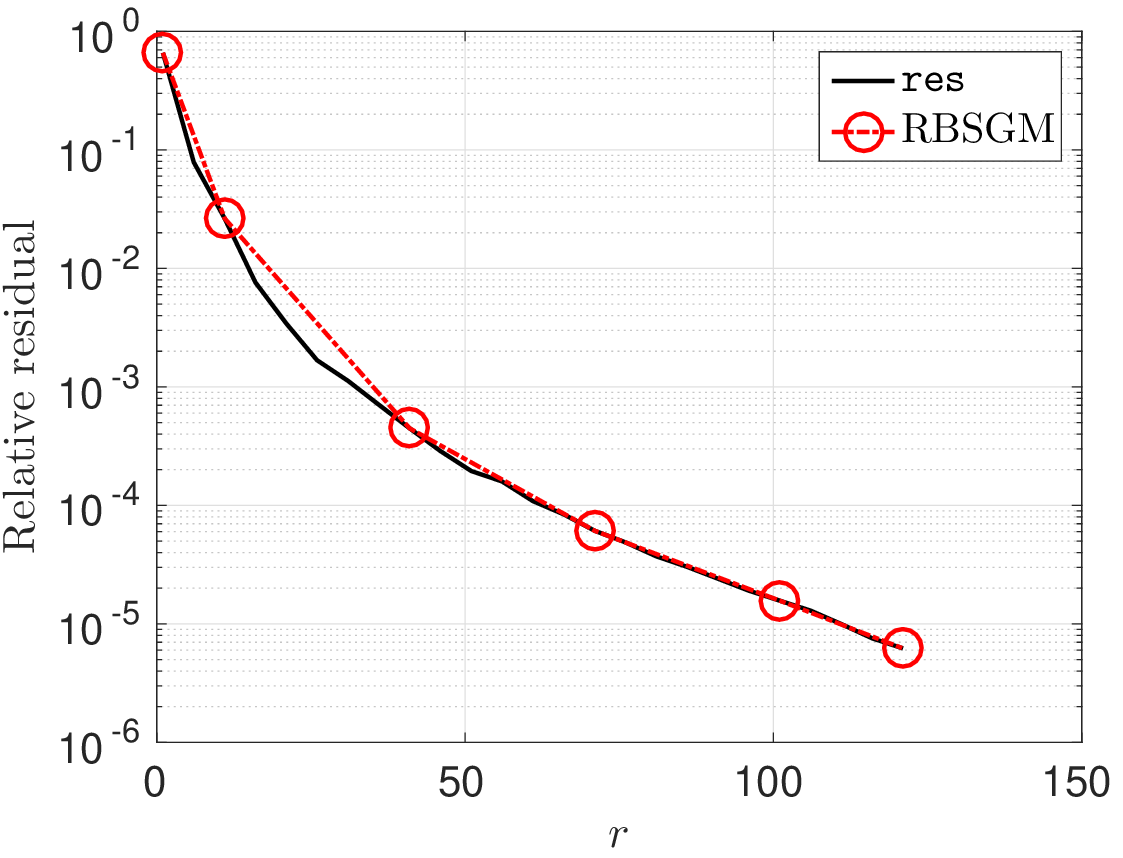}
		}
		\subfloat[$m=7,N_h=129^2$]{
			\includegraphics[width=0.32\linewidth]{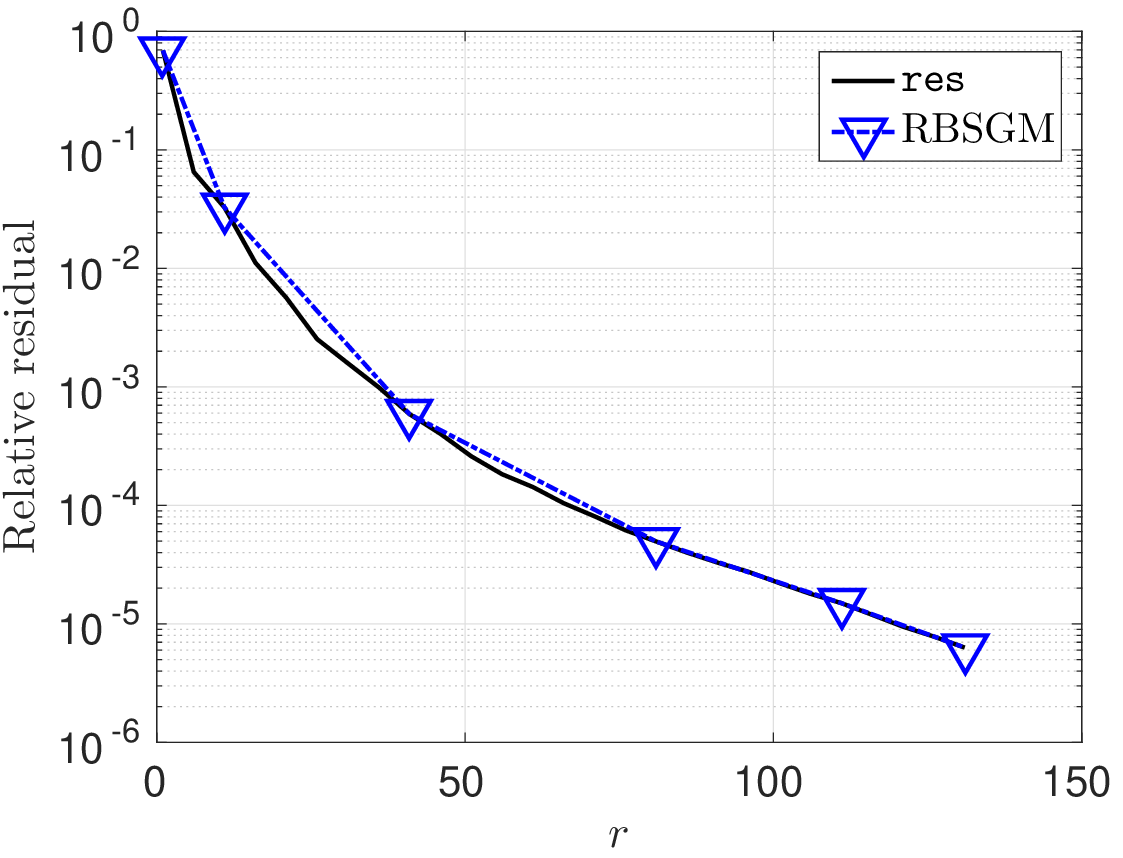}
		}\\
		\subfloat[$m=10,N_h=33^2$]{\includegraphics[width=0.32\linewidth]{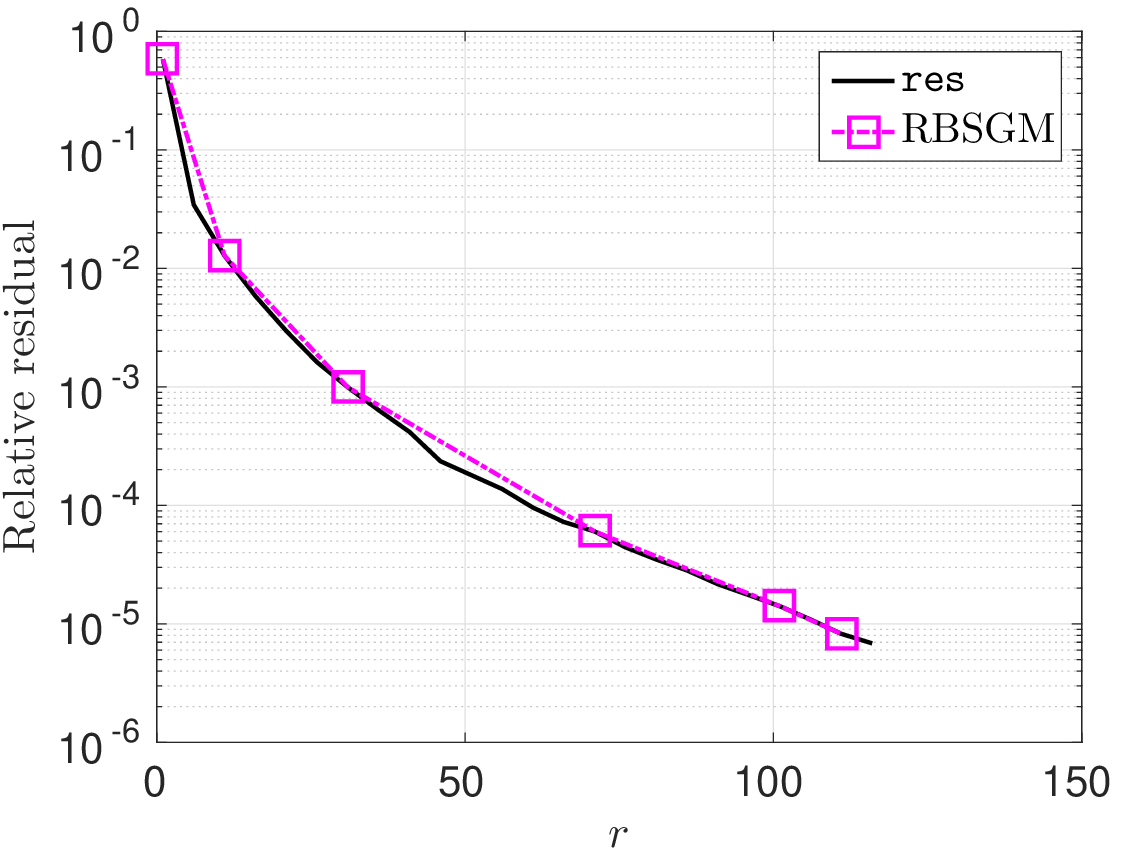}
		}
		\subfloat[$m=10,N_h=65^2$]{
			\includegraphics[width=0.32\linewidth]{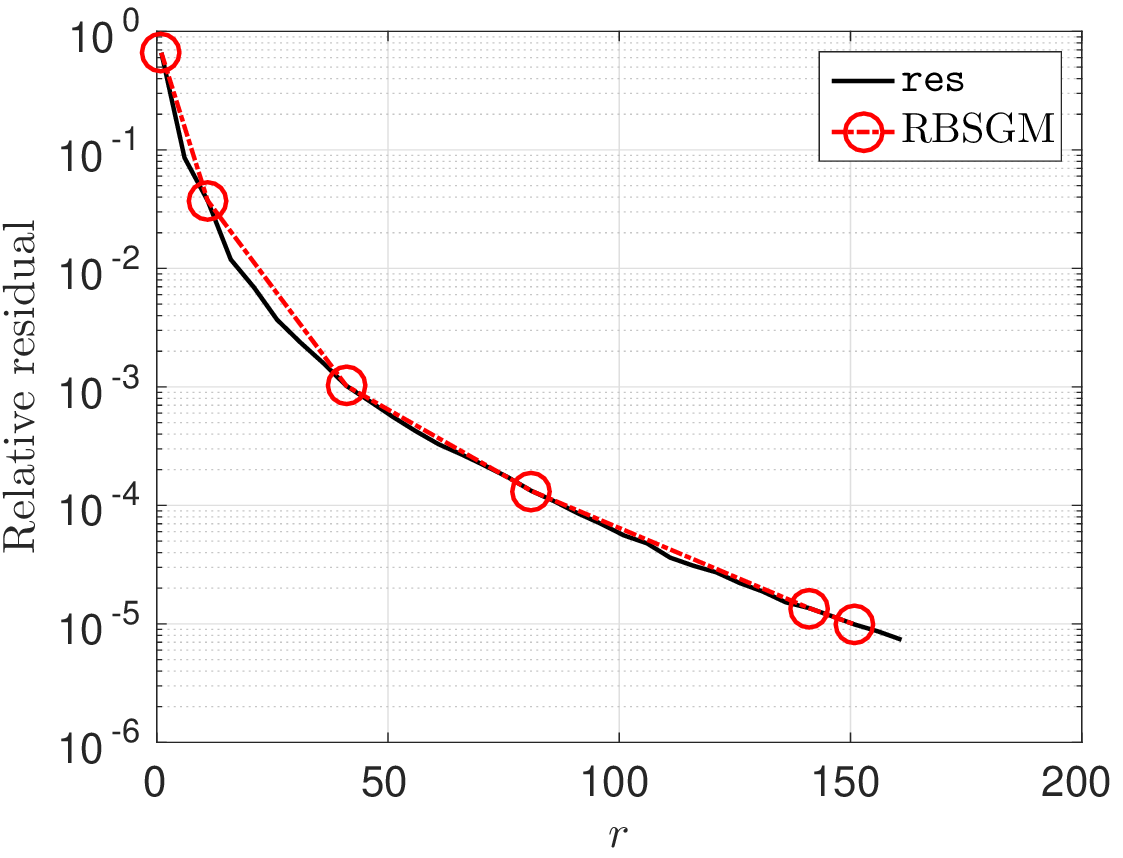}
		}
		\subfloat[$m=10,N_h=129^2$]{
			\includegraphics[width=0.32\linewidth]{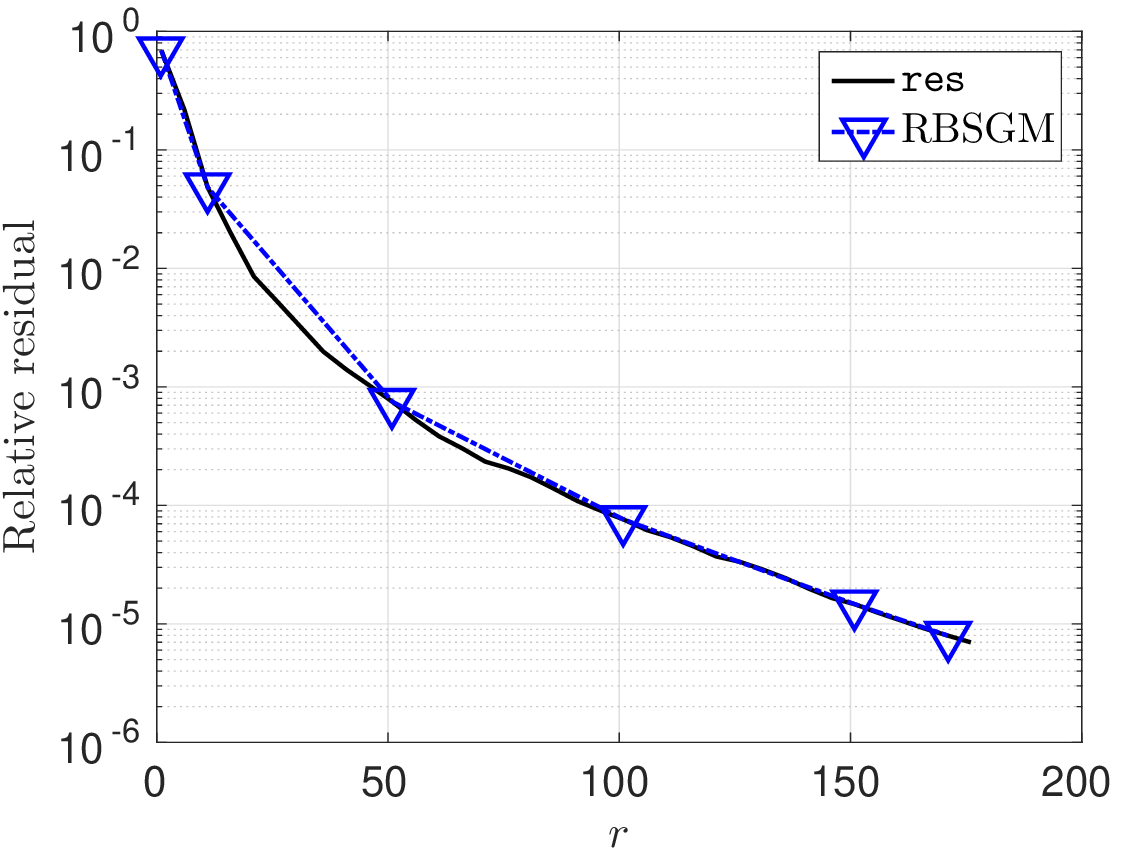}
		}\\
	\end{center}
	\caption{Relative residual with respect to the number of reduced basis functions $r$.}\label{fig:test2a}
\end{figure}
In Figure~\ref{fig:test2a}, we present the relationship between the number of reduced basis functions and the corresponding relative residual. Additionally, we provide the relative residual computed during the adaptive procedure of the RBSGM. Both the full linear system~\eqref{eq:linear_sg} and the reduced linear system~\eqref{eq:rbsgeq} are solved using the Bi-CGSTAB with a mean based preconditioner. The figure illustrates that the number of reduced basis functions, denoted by $r$, is significantly smaller than the DOF in the finite element method. Furthermore, it is observed that the value of $r$ is influenced by the chosen tolerance $\tol$ but shows little dependence on the grid resolution $\Nx$. In the RBSGM, the relative residual only needs to be computed $5$ or $6$ times when using $\tol = 10^{-5}$. This emphasizes the computational efficiency of the reduced basis stochastic Galerkin method, as it achieves comparable accuracy with significantly fewer computational costs.

\begin{table}[htbp]
	\caption{Relative errors of the mean and variance functions.}\label{tab:test2c}
\begin{center}
	\newcolumntype{C}{>{\centering\arraybackslash}X}%
	\begin{tabularx}{0.95\linewidth}{c|c|cCCC}
		\hline
		$\tol$ & & $\Nx$ & $m = 5$  & $m=7$ & $m=10$\\
		\hline
		\multirow{6}*{$10^{-4}$}&	& $33^2$& 4.12e-06 (2.53e-06)  & 7.73e-06 (2.96e-06)  & 7.96e-06 (3.11e-06)  \\
		& $\meanerr$& $65^2$& 1.22e-05 (4.86e-06)  & 5.06e-06 (5.21e-06)  & 6.52e-06 (5.40e-06)  \\  
		& & $129^2$& 7.21e-06 (1.30e-05)  & 9.22e-06 (1.33e-05)  & 6.22e-06 (1.40e-05)  \\ 
		\cline{2-6} 
		& & $33^2$& 2.04e-05 (4.01e-05)  & 3.32e-05 (4.41e-05)  & 3.05e-05 (4.56e-05)  \\
		& $\varerr$& $65^2$& 4.62e-05 (1.30e-05)  & 1.77e-05 (1.42e-05)  & 2.02e-05 (1.61e-05)  \\  
		& & $129^2$& 2.09e-05 (2.16e-05)  & 2.25e-05 (2.14e-05)  & 1.63e-05 (2.28e-05)  \\  
		\hline 
		\multirow{6}*{$10^{-5}$}		& &$33$& 6.84e-07 (1.32e-06)  & 4.22e-07 (1.50e-06)  & 4.25e-07 (1.53e-06)  \\  
		&$\meanerr$ &$65$& 8.01e-07 (7.25e-07)  & 2.87e-07 (9.95e-07)  & 5.69e-07 (2.52e-07)  \\  
		& &$129$& 5.04e-07 (7.71e-07)  & 4.15e-07 (8.35e-07)  & 3.62e-07 (1.00e-06)  \\
		\cline{2-6}
		& &$33^2$& 2.83e-06 (4.17e-06)  & 1.67e-06 (4.75e-06)  & 1.24e-06 (4.98e-06)  \\  
		&$\varerr$ &$65^2$& 2.15e-06 (2.99e-06)  & 9.43e-07 (4.30e-06)  & 1.96e-06 (1.50e-06)  \\  
		& &$129^2$& 1.90e-06 (2.02e-06)  & 1.08e-06 (1.47e-06)  & 1.10e-06 (1.60e-06)  \\  
		\hline
	\end{tabularx}
\end{center}
\end{table}

Table~\ref{tab:test2c} presents the relative errors of the mean and variance functions for RBSGM. In addition, we provide the relative errors of the SGM in parentheses for comparison. The table reveals that the relative errors of the mean and variance functions in RBSGM are comparable to those in the SGM.  It is worth noting that, in the adaptive procedure, we solve the reduced basis linear system~\eqref{eq:rbsgeq} using the Bi-CGSTAB with a tolerance of  $10^{-7}$. This tolerance value is smaller than $\tol$ used in this test problem. Meanwhile, since we increase the number of reduced basis functions by a multiple of $\mathtt{ns}=10$ during each adaptive iteration in Algorithm~\ref{alg:RBGM}, the number of reduced basis functions determined by the adaptive procedure may be slightly larger than that required. As a result, RBSGM outperforms SGM in terms of accuracy for many cases in this test problem.

\section{Conclusions}
In this work, we develop a reduced basis stochastic Galerkin method for partial differential equations with random inputs. In comparison to the standard stochastic Galerkin method, which discretizes the physical space using grid based approaches, our proposed method offers improved computational efficiency in computing the Galerkin solution, particularly when dealing with a large number of physical degrees of freedom. While our current focus has been on efficiently identifying reduced basis functions for the physical approximation, our future work will concentrate on constructing effective bases for the stochastic approximation.

\section*{Acknowledgments}
This work is supported by the National Natural Science Foundation of China (No. 12071291), the Science and Technology Commission of Shanghai Municipality (No. 20JC1414300) and the Natural Science Foundation of Shanghai (No. 20ZR1436200)


  \bibliographystyle{elsarticle-num} 
  \bibliography{ref}



%

%
%
%
\end{document}